\documentclass[11pt]{mcom-l}

\usepackage[titletoc]{appendix}

\usepackage{amsmath}
\usepackage{amsfonts}
\usepackage{color}
\usepackage{amsrefs}
\usepackage{multirow}
\usepackage{graphicx}
\usepackage{enumitem}
\usepackage{amssymb}
\usepackage[normalem]{ulem}

\newtheorem{theorem}{Theorem}[section]    
\newtheorem{cor}{Corollary}[section]    
\newtheorem{lemma}{Lemma}[section]

\newtheorem{assump}{Assumption}

\newcommand{\R}{\mathbb{R}}

\providecommand{\norm}[1]{\lVert#1 \rVert}

\newcommand{\floor}[1]{\left\lfloor #1 \right\rfloor}

\numberwithin{equation}{section}
\newcommand{\halmos}{\hfill$\Box$}

\newcommand{\yrest}{y_R^k}

\begin{document}

\title[Inexact Restoration for Inexact Problems]{Inexact Restoration for Minimization with Inexact Evaluation both of the Objective Function and the Constraints}

\author[]{L. F. Bueno}
\address{Institute of Science and Technology, Federal University of S\~ao Paulo, S\~ao Jos\'e dos Campos SP, Brazil. e-mail:  lfelipebueno@gmail.com }
\author[]{F. Larreal}
\address{{Department of Applied Mathematics, Institute of Mathematics, Statistics, and Scientific Computing (IMECC), State University of Campinas, 13083-859 Campinas SP, Brazil. e-mail: francislarreal@gmail.com}}
\author[]{J. M. Mart\'{\i}nez}
\address{Department of Applied Mathematics, Institute of Mathematics, Statistics, and Scientific Computing (IMECC), State University of Campinas, 13083-859 Campinas SP, Brazil. e-mail: martinez@ime.unicamp.br}
\thanks{This work was supported by FAPESP (grants
    2013/07375-0, 2018/24293-0 and 2021/14011-0), CAPES and CNPq.}
\date{January 02, 2022, Revised September 19, 2023}
\maketitle

\begin{abstract}
 In a recent paper an Inexact Restoration method for solving 
 continuous constrained optimization problems was analyzed from the point
  of view of worst-case functional complexity and convergence. On the other hand, 
  the Inexact Restoration methodology was employed, in a different research,
     to handle minimization problems
   with inexact evaluation and simple constraints. These two methodologies are 
   combined in the present report, for constrained minimization problems in which
   both the objective function and the constraints, as well as their derivatives, 
    are subject to evaluation errors. 
  Together with a complete description of the method, complexity and convergence 
   results will be proved. 
  
\end{abstract}

\section{Introduction} \label{introduction}

 Consider an optimization problem given by

     \begin{equation} \label{theproblem}
        \begin{array}{cl} 
        \mbox{Minimize }   &  F(x) \\ 
        \mbox{ subject to }  &   H(x)= 0 \\
                            & x \in \Omega,
        \end{array}
    \end{equation}
   where  $ F: \R^{n} \to \R $, $ H: \R^n \to \R^m $ and
   $ \Omega $ is a nonempty compact polytope. As usually, if inequality constraints
  $G(x) \leq 0$ are present, we reduce the problem to the standard form (\ref{theproblem}) by means 
  of the addition of slack variables. Assume that exact evaluation of 
  $F(x)$, $H(x)$ and their derivatives is not always possible. Instead, each evaluation 
   of $F(x)$ (or $H(x)$) is, according to availability or convenience, replaced with 
  $f(x, y)$ (or $h(x, y)$, respectively) where $y$ lies in an abstract set $Y$ and 
  determines the degree of precision in which the objective function or the constraints
  are evaluated. We will assume that $g_f: Y \to \R_+$ is such that $f(x, y) = F(x)$ 
  when $g_f(y)=0$, $g_h: Y \to \R_+$ is such that $h(x, y) = H(x)$ when $g_h(y)=0$   
   and that, roughly speaking, the precision in the 
   evaluations improves when $g_f(y)$ and $g_h(y)$ decrease. If 
 the precision of the objective function is governed by a set $Y_1$   and the precision of the constraints are
 governed by $Y_2$, where both $Y_1$ and $Y_2$ are abstract sets, we may define  
 $Y = Y_1 \times Y_2$.  
  Writing $g(y)=\max\{g_f(y),g_h(y)\}$, problem (\ref{theproblem}) is equivalent to:
\begin{equation} \label{eq:problemainexato}
\begin{array}{cl} 
\mbox{Minimize (with respect to } x \mbox{)}&  f(x,y)  \\ 
\mbox{ subject to }&  h(x,y) = 0, \\
& g(y) = 0, \\
& x \in \Omega,\\
& y \in Y.
\end{array}
\end{equation} 

  A solution of (\ref{eq:problemainexato}) could be obtained fixing $y \in \Omega$ in such 
   a way 
  that $g_f(y) = g_h(y) = 0$ and handling the resulting problem as a standard constrained
  optimization problem. However, we are interested in problems in which such procedure 
   is not affordable because solving (\ref{eq:problemainexato}) fixing $g_f(y)=g_h(y)=0$ 
  is overwhelmingly expensive or even impossible.    

  The definition (\ref{eq:problemainexato}) makes sense independently of the meaning of 
  $y, Y$, or $g(y)$. We have especially in mind the case in which $f(x, y)$ 
   represents $F(x)$ with an error governed by $y \in Y$, $h(x, y)$ is $H(x)$ computed with
  an error that depends on $y$, $g_f(y)=0$ if and only if $f(x, y) = F(x)$, and 
   $g_h(y)=0$ if and only if 
  $h(x, y) = H(x)$ for all $x \in \Omega$. However, the results of this paper can be read
  without reference to this meaning. 

 In this paper we extend the results of \cite{bkm2019} and \cite{ir22}. In \cite{bkm2019} the
  problem (\ref{eq:problemainexato}) is considered without the presence of the constraints
  $h(x, y)=0$. In \cite{ir22} an Inexact Restoration method with worst-case complexity
  results is introduced for solving the classical constrained optimization problem. The 
   techniques of \cite{bkm2019} and \cite{ir22} are merged in the present paper in order 
   to handle the constrained optimization problem with inexactness both in the objective function
  and the constraints.

 Let us give an example of the applicability of the present approach which, in fact,
  motivated the algorithmic framework and theoretical analysis developed in this paper.  We are involved with real-life river simulations
and the corresponding inverse problems \cite{livrocriab} and \cite{saintvenant}.        
 The Saint-Venant  equations 
\begin{equation} \label{massa}
\frac{\partial A}{\partial t} + \frac{\partial Q}{\partial x} = 0
\end{equation}
and
\begin{equation} \label{momentong}
\frac{\partial Q}{\partial t} + \frac{\partial}{\partial x}
\left(\frac{Q^2}{A}\right) + {\tt g} A \frac{\partial z}{\partial x} +
\frac{{\tt n}^2 g Q|Q|}{A R^{4/3}} = 0
\end{equation} 
 are usually employed for
 river-flow simulations. In (\ref{massa}) and (\ref{momentong})       
 $x \in [x_{\min},x_{\max}]$ and $t \in [t_{\min},t_{\max}]$, where  $z_b(x)$ is the 
bed elevation, measured from a
  datum, 
$z(x,t) - z_b(x)$ is the depth of the river at $(x,t)$,
$A(x,t) = [z(x, t)-z_b(x, t)], w(x)$ is the transversal wetted area at $(x,t)$,
$P(x,t) = w(x) + 2 [z(x, t)-z_b(x, t)]$ is the wetted perimeter at $(x,t)$, $R(x,t)
= A(x,t) / P(x,t)$ is the hydraulics radius at $(x,t)$, $V(x,t) =
Q(x,t) / A(x,t)$ is the speed of the fluid at $(x,t)$, and ${\tt g}$ is the
acceleration of gravity.  Equation~(\ref{massa})
describes mass conservation and equation~(\ref{momentong}) represents
conservation of the linear momentum. 
      Finally, 
 ${\tt  n}$ is called Manning Coefficient and takes account of friction.

 When the Saint-Venant equations are solved by means of a stable implicit method \cite{leveque}, the 
estimation of Manning coefficients require 
to solve a constrained optimization problem \cite{ayvaz}.

 However real rivers are not rectilinear, their flux is not homogeneous, cross sections are not rectangular and, sometimes, 
  are time-dependent, and the Manning coefficients are not constant. Therefore, increasing levels of problem complexity
  arise when we incorporate variable Manning coefficients with different dimensions, cross section variations and when 
  we increase the number of observations or expert guesses for the flux evolution. Further difficulties arise when we 
  refine the discretization grid for solving realistic Saint Venant equations. These considerations lead to different 
   formulations of real-life river simulations, each corresponding to variable precisions for the computation of
   the objective function and the constraints.

  Inexact Restoration (IR) methods for constrained continuous optimization
     were introduced in \cite{ir2}, 
  inspired  
   in several classical papers by Rosen \cite{rosen} and Miele \cite{mhh}, among others. 
  At each iteration of an IR algorithm feasibility is firstly improved and, then, optimality
   is improved along a tangent approximation of the feasible region. The 
   so far generated trial point is accepted or not as new iterate according to the decrease
   of a merit function or using filter criteria \cites{gkv,ir6,ir17,ir18,ir21}.
 Theoretical papers concerning
   Inexact Restoration methods for constrained optimization include \cites{ir4,ir15,ir22}.
    Algorithmic variations are discussed in \cites{ir4,bfms,ir3,ir1,ir16,ir17,ir18,ir25}, and applications
  may be found in \cites{ir5,ir7,ir19,ir8,ir9,ir10,ir13,ir12,ir14,ir20,ir27,ir24,ir26}. 

   The idea of using the IR framework to deal with optimization problems in which 
  the objective function is subject to evaluation errors comes from \cite{nkjmm}, where 
  inexactness came from the fact that the evaluation was the result of an iterative process. 
    Evaluating the function with additional precision was considered in \cite{nkjmm} as a
  sort of inexact restoration. This basic principle was developed in \cite{bkm2018} and 
  \cite{bkm2019}, where complexity results were also proved. Moreover, in 
  \cite{bkm2021} the case in which derivatives are not available was considered.
     Inexactness of the objective function in optimization problems was addressed in 
     several additional papers in recent years \cites{inexact1,inexact2,ir23,inexact7,inexact4,inexact6,inexact3}.   
   The objective of the present paper is to use the ideas of \cites{bkm2018,bkm2019,nkjmm} to 
   handle the constrained optimization problem in which the evaluation of the objective 
   functions and the constraints is subject to error. We will show that, although the 
   main ideas are applicable, a number of technical difficulties appear whose solution
   offer additional insight in the problem. From the theoretical point of view we will prove
   convergence to feasible points (whenever possible) and asymptotic fulfillment of 
   optimality conditions.  

 This paper is organized as follows. In Section~\ref{algorithms} we describe BIRA, 
    the main algorithm for solving \eqref{eq:problemainexato}. In Section~\ref{purpose} we 
  state our final goal in terms of complexity and convergence of BIRA and we highlight the general
 lines that will be followed in the proofs. In Section~\ref{assumptions} 
   we state general assumptions on the problem that will be used throughout the paper.   
    In Section~\ref{restauracao} we state several theoretical results with respect to the
  Restoration Algorithm RESTA that will be useful in forthcoming sections. In Section~\ref{birawell} 
   we show that every iteration of BIRA is well defined. In Section~\ref{conviastep} we 
  prove convergence towards feasible points. 
    In Section~\ref{compcon} we finish up proving complexity and convergence
  of the main algorithm. Conclusions are stated in
 Section~\ref{conclusions}. Proofs of the technical lemmas are presented in Appendix A. \\

\noindent
{\bf Notation }   

\begin{enumerate}

\item  All along this paper $\|\cdot\|$ represents the Euclidean norm.

\item We define 
\begin{equation} \label{defc}
c(x, y) = \frac{1}{2}\|h(x, y)\|^2.\\      
\end{equation}    
 
\item $P_\Omega(z)$ denotes the Euclidean projection of $z$ onto $\Omega$. 

\end{enumerate}

\section{Algorithms}  \label{algorithms}

  In this section we define the Basic Inexact Restoration Algorithm (BIRA) for solving our main
  problem and the restoration algorithm RESTA, which is called at each iteration
   of BIRA.
  
All along the paper we will use the merit function that combines objective function and constraints
  defined by means of the penalty parameter $\theta$ according to:
     \begin{equation} \label{defmeritfunction}
         \Phi(x,y, \theta) = \theta f(x,y)  + (1 - \theta)\left[ \|h(x,y)\| +g(y)\right]
    \end{equation}
  for all $x \in \Omega, y \in Y$, and $\theta \in [0, 1]$.

     \subsection{Basic Inexact Restoration Algorithm (BIRA)}

    The iterative Algorithm BIRA has  three main steps. Each iteration begins with a Restoration
  Phase, at which, starting from the current iteration $x^k$ and the 
  current precision variable $y^k$, one computes an inexactly restored
   $x_R^k$ and a better precision parameter $y^k_R$. At the second step, the penalty parameter
   that defines the merit function is conveniently updated.
 At Step~3 (Optimization Phase) we try to improve the merit function 
    by approximate minimization of a quadratic approximation of the objective function with
  an adaptive regularization parameter.  At the first iterations of the Optimization Phase 
   we admit to relax the accuracy defined by $y^k_R$ with the aim of reducing computational cost. If this 
   relaxation is not successful the Optimization Phase uses the precision level $y^k_R$.

   The description of Algorithm BIRA begins reporting all the algorithmic parameters that will
   be used in the calculations. The parameter $r \in (0, 1)$ is used in the Restoration Phase.
   At Step~2 we use $\theta_0 \in (0, 1)$, the initial penalty parameter. Bounds for the first
  regularization parameter used in the Optimization Phase are given by $\mu_{min}$ and 
  $\mu_{max}$.The parameter 
  $\alpha > 0$ is used at Step~3 to decide acceptance or rejection of the trial point 
    obtained at this step, $M$ is a bound for Hessian approximations, and $N_{acce}$ is the 
    maximal number of steps, at the Optimization Phase, in which relaxing precision is admitted.  

   Other parameters ($\alpha_R, \sigma_{max}, \sigma_{min},  \beta_c, r_{feas}, 
  \bar{\epsilon}_{prec}, N_{prec}, \beta_{PDP}$) are used in Algorithm RESTA  and 
   will be commented later.\\

\noindent
\textbf{ Algorithm \ref{algorithms}.1 (BIRA)}

Given  $ \alpha_R, \alpha> 0 $, $ M \geq 1 $, $ \sigma_{max} \geq \sigma_{min}> 0 $,
 $ \mu_{max}\geq \mu_{min}> 0 $, $\beta_c>0$, $\beta_{PDP}>0$, 
     $ r \in (0, 1) $, $ r_{feas} \in (0, r) $, $\bar{\epsilon}_{prec} \geq 0$, $N_{prec} \geq 0$, and  
  $N_{acce} \geq 0$, 
   choose  $\mu_{-1} \in  [\mu_{min}, \mu_{max}]$,
      $ x^0 \in \Omega$, $y^0 \in Y $, set $ k \leftarrow  0 $ and  
 $ \theta_0 \in (0,1) $.

\noindent
{\bf Step 1. Restoration Phase}

   Compute 
 $(x_R^k,y_R^k)$ using Algorithm RESTA. 
 
  If 
\begin{equation} \label{fallar}
 \|h(x_R^k,y_R^k)\|  > r \, \|h(x^k,y_R^k)\|
\end{equation}
or
	\begin{equation}\label{errata}
		\|h(x^k, y_R^k)\|-\|h(x_R^k, y_R^k)\| < \dfrac{1-r}{2r}\left[  g(y^k)-g(y_R^k) \right],	
\end{equation}
 stop the execution of BIRA declaring
  \textit{Restoration Failure}.

\noindent
{\bf Step 2. Update penalty parameter}

Test the inequality 
        \begin{equation} \label{checkpenalty}
        \begin{array}{c}
        \Phi(x_R^k, y_R^k,   \theta_k)-\Phi(x^k, y^k_R,   \theta_k) \\
        \le   \frac{1-r}{2}\left[  \|h(x_R^k,y_R^k)\|   - \|h(x^k,y^k_R)\|+  g(y_R^k)-g(y^k)\right].
        \end{array}            
        \end{equation} 
     If (\ref{checkpenalty})    holds, define $\theta_{k+1} = \theta_k$. 
         Else, compute  
        \begin{equation} \label{updatepenaltyparameter}
        \begin{array}{c}
 {\theta}_{k+1} = 
 \frac{(1+r) \left[  \|h(x^k,y^k_R)\|  - \|h(x_R^k,y_R^k)\| + g(y^k) - g(y_R^k)\right]} { 2 \left[ f(x_R^k,y_R^k) - f(x^k,y^k_R) +\|h(x^k,y^k_R)\| -  \|h(x_R^k,y_R^k)\|  + g(y^k) - g(y_R^k)\right]}.
        \end{array} 
        \end{equation}

\noindent
{\bf Step 3. Optimization Phase} 

   Initialize $\ell \leftarrow 0$.

\noindent
{\bf Step 3.1} Choose $y^{k+1} \in Y$ (perhaps $g(y^{k+1})$ bigger than  $g(y_R^k)$). Choose
   $\mu\in  [\mu_{min}, \mu_{max}]$ and  a symmetric matrix $H_k \in \mathbb{R}^{n \times n}$ such that $\norm{H_k} \leq M$. 

\noindent
{\bf Step 3.2}    
      If  $\ell \geq N_{acce}$, re-define $y^{k+1} = y^k_R$.

\noindent
{\bf Step 3.3 } Compute 
	 $x \in \Omega$  an approximate solution of        
	\begin{equation} \label{subproblemintermediario}
		\begin{array}{ll}
			\mbox{Min.} & \nabla_x f(x_R^k,y^{k+1})^T (x - x_R^k)  + \frac{1}{2} (x- x_R^k)^T H_k (x - x_R^k) + \mu \norm{ x - x_R^k}^2 \\
			\mbox{s. to }  &\nabla_x h(x_R^k, y^{k+1})^T(x-x_R^k) = 0\\     & x \in \Omega.  
		\end{array}
	\end{equation} 

\noindent
{\bf Step 3.4  } Test the conditions
	\begin{equation}\label{desfinte}
	f(x,y^{k+1}) \leq f(x_R^k,y_R^k) - \alpha \norm{x - x_R^k}^2
	\end{equation}		
	 and 
			 \begin{equation}\label{desmeri}
    \begin{split}
        \Phi(x, y^{k+1}, \theta_{k+1})  
   &\leq  \Phi(x^k, y^{k+1}, \theta_{k+1}) \\& \quad + \frac{1-r}{2}\left[ \|h(x^k_R,y_R^k)\|  - \|h(x^k,y^k_R)\|  + g(y_R^k)- g(y^k) \right]. 
    \end{split}
		\end{equation} 
		
	     If   \eqref{desfinte} and \eqref{desmeri} are fulfilled, 
   define $\mu_k = \mu$, $x^{k+1} = x$, $k \leftarrow k+1$, and go to Step~1. 
		Otherwise, choose  $\mu_{new} \in [2 \mu, 10 \mu]$, $\mu \leftarrow \mu_{new}$,
 set  $\ell \leftarrow \ell+1$ and go to Step 3.2.\\

\noindent
{\bf Remark } 
In 
   Assumption~\ref{assumpprecsubprob}     
               we will define in which sense, at Step~3.3, $x$ should be an approximate 
  solution of (\ref{subproblemintermediario}). 

The condition (\ref{fallar}) used at Step~1 in BIRA is not the natural generalization 
  of the restoration condition used in previous IR algorithms. Such ``natural" generalization should be
  $\|h(x^k_R, y^k_R)\| > r \|h(x^k, y^k)\|$. The reason why the traditional alternative is not adequate in the context
  of BIRA is the following: Suppose that, by chance, 
  $\|h(x^k, y^k)\|$ vanishes or is very small. In this case, the restored $(x^k_R, y^k_R)$ would be rejected almost certainly, and
   the algorithm would stop by Restoration Failure. However, this decision could be unreasonable because, even if 
  $\|h(x^k_R, y^k_R)\|$ is greater than $\|h(x^k, y^k)\|$, the point $x^k_R$ could be better than 
   $x^k$  when the constraints are evaluated with the same accuracy defined by $y^k_R$, which may be substantially better than the one defined by $y^k$. 
    This is the reason why we preferred (\ref{fallar}) instead of   $\|h(x^k_R, y^k_R)\| > r \|h(x^k, y^k)\|$ for 
   deciding failure of restoration. In general, the level of precision used in each of the conditions used in the algorithm needs to be carefully chosen. A technical consequence of these decisions is that the theoretical proofs 
    in this paper are, many times, reasonably  different than the corresponding proofs of other IR papers.

\subsection{Algorithm for the Restoration Phase} 

  The objective of the restoration algorithm RESTA is to find $x^k_R$ and $y^k_R$ such that 
 the inequalities (\ref{erre3}) below are fulfilled. 
  In general, the fulfillment of  $   g_f(y^k_R) \leq r g_f(y^k)$ and 
    $g_h(y^k_R) \leq r g_h(y^k)$  is easy to obtain as, under the usual interpretation,  
  these inequalities merely impose that the precision with which $F$ and $H$ are evaluated at 
  $x^k_R$ should be better than the precision with which $F$ and $H$ were evaluated at $x^k$. 
  However, the requirement  $\|h(x_R^k,y_R^k)\|  \leq r \, \|h(x^k,y_R^k)\|$  could be difficult to achieve. We try to do this by minimizing a regularized quadratic model of the sum of squares of the constraints. The regularization parameter is initialized between $\sigma_{min}$ and $\sigma_{max}$ and $\alpha_R$ is associated with the sufficient decrease criterion for acceptance of the trial point. Parameters $\beta_c$ and $\beta_{PDP}$ control the distance between some restored point estimates and the current iterate.

  In critical cases,
  where the original
   problem is infeasible, the fulfillment of  $\|h(x_R^k,y_R^k)\|  \leq r \, \|h(x^k,y_R^k)\|$  could be even impossible. Therefore, 
   ``Restoration Failure" is a possible diagnostic that needs to appear  
   at any algorithm that aims to fulfill those requirements. In order to declare that we are probably in this situation, we use the parameter $r_{feas}$, defined to be smaller than $r$ in BIRA, to check if the projected gradient of the sum of squares of the constraint violations is sufficiently smaller than the infeasibility measure. When we solve the problem with
   precision $w^i$ and we obtain a point $z^l$ indicating that the original problem may be infeasible, 
   we have to decide whether we progressively try to get out of this situation by improving precision and 
   seeking a smaller infeasibility with respect to $h$ or if we demand more quickly a better quality in 
   the representation of constraints and their derivatives, decreasing $g_h$, to accurately check the infeasibility status. 
   The $N_{prec}$ parameter  determines a limit of attempts with an indication
   of infeasibility until we force the precision in the calculation of the constraints to be at $\bar{\epsilon}_{prec}$,  the level required by the user.  \\

 \noindent
\textbf{ Algorithm \ref{algorithms}.2 (RESTA)}

  Assume that  $x^k \in \Omega$, $y^k \in Y$, and the parameters that define BIRA are
  given. If 
     $\|h(x^k,y^k)\| + g(y^k) =0 $, return defining $(x_R^k,y_R^k)=(x^k,y^k)$. Else, set $i \leftarrow 0$ and $w^0 = y^k$.   
    
\noindent
 {\bf Step 1 } 
  Using an optional inexpensive
  problem-dependent procedure (PDP) (if available), try to compute $y_R^k \in Y$ and 
  $x_R^k \in \Omega$ such that 
\begin{equation} \label{erre3}
 g_f(y^k_R) \leq r g_f(y^k), \;   g_h(y^k_R) \leq r g_h(y^k), \;  
   \|h(x_R^k,y_R^k)\|  \leq r \|h(x^k,y_R^k)\|,    
\end{equation}
and
\begin{equation} \label{beta3}
\max\{  \norm{ x_R^k-x^k}, \norm{ y_R^k-y^k} \}\leq \beta_{PDP}  \|h(x^k,y^k_R)\| . 
\end{equation} 
  If such procedure is activated and (\ref{erre3}) and (\ref{beta3}) hold, return.

\noindent
{\bf Step 2 }
    	  If $i \leq N_{prec}$, set $\bar{g}_h \leftarrow r g_h(w^i)$, else
  $\bar{g}_h \leftarrow \min\{\bar{\epsilon}_{prec},r g_h(w^i)\}$.  
   If  $g_f(w^i)=0$ and $g_h(w^i)=0$ define $w^{i+1}=w^i$, else choose  $w^{i+1} \in Y$ such that 
    	\begin{equation}\label{restorationg} 
    	g_f(w^{i+1}) \leq r g_f(y^k) \,\,\,\, \text{ and } \,\,\, g_h(w^{i+1}) \leq \bar{g}_h.
    	\end{equation} 

(This choice of $w^{i+1}$ will be assumed to be possible and inexpensive since, in general, 
  merely represents increasing the precision of forthcoming evaluations.)

\noindent
{\bf Step 3 } Compute   
    	 $z^0 \in \Omega$ such that
\begin{equation} \label{cz1}
 c(z^0,w^{i+1}) \leq c(x^k,w^{i+1})
\end{equation} 
(see (\ref{defc}) for the definition of $c$)  and
\begin{equation} \label{cz2}
  \|z^0-x^k\| \leq \beta_c \|h(x^k,w^{i+1})\|.
\end{equation}
 (Note that the choice of $z^0$ satisfying (\ref{cz1}) and (\ref{cz2}) is always possible because
  the trivial choice $z^0 = x^k$ is admissible.)

     Set $\ell \leftarrow  0$.

\noindent
{\bf Step 4 } Test the stopping criteria 
\begin{equation} \label{thesto}
  c(z^\ell,w^{i+1}) \leq  r^2  c(x^k,w^{i+1})
\end{equation}
and
 \begin{equation} \label{thesto2}
    	\|P_{\Omega} \left (z^\ell-\nabla_x c(z^\ell,w^{i+1}\right)-z^\ell \| \leq 
    r_{feas} \|h(x^k,w^{i+1})\| 
  \mbox{ and } g_h(w^{i+1}) \leq \bar{\epsilon}_{prec}.
\end{equation}

If (\ref{thesto}) holds or (\ref{thesto2}) holds, return to BIRA defining   
    	 $x_R^k=z^\ell$ and $y^k_R=w^{i+1}$.

 (Although both (\ref{thesto}) and (\ref{thesto2}) are reasons for returning, 
  these inequalities have quite different meanings since (\ref{thesto}) indicates 
   success of restoration whereas (\ref{thesto2}) indicates possible failure. 
  In any case, the final success restoration test is made in BIRA.)       

If 
\begin{equation} \label{thesto3}
 \|P_{\Omega} \left (z^\ell-\nabla_x c(z^\ell,w^{i+1}\right)-z^\ell \| 
\leq \epsilon_c \mbox{ and } g_h(w^{i+1}) > \bar{\epsilon}_{prec},
\end{equation}
 set $i \leftarrow i+1$ and go to Step~2.

\noindent
{\bf Step 5} Choose 
     $\sigma \in  [\sigma_{min}, \sigma_{max}]$ and $B_\ell \in \mathbb{R}^{n \times n}$    
  such that $B_\ell+ \sigma_{min}  I$ be
 symmetric and positive definite with  $\|B_\ell \|\leq M$ and  $ \|{(B_\ell+ \sigma_{min} I)^{-1}}\|  \leq M$.

\noindent
{\bf Step 5.1} 
 Compute 
	 $z^{trial} \in \Omega$  as an approximate solution of 
    		\begin{equation} \label{subproblemfeas}
    		\begin{array}{ll}
    		\text{ Minimize} & \nabla_x c(z^\ell,w^{i+1})^T (z - z^\ell) + \frac{1}{2} (z- z^\ell)^T (B_\ell+ \sigma I) (z- z^\ell)  \\
    		\text{ subject to } & z \in \Omega.
    		\end{array}
    		\end{equation}

\noindent
{\bf Step 5.2}  
 Test the condition
    		\begin{equation} \label{desc}
    		c(z^{trial},w^{i+1}) \leq c(z^\ell,w^{i+1}) - \alpha_R \|{z^{trial}- z^\ell}\|^2.
    		\end{equation}
    		If  \eqref{desc} is fulfilled, define $z^{\ell+1}  = z^{trial}$,   
    	      set  $\ell \leftarrow \ell+1$,  and go to Step~4. Otherwise, choose            
    		\begin{equation} \label{updatemufeas}
  		 \sigma_{new}  \in [2 \sigma, 10 \sigma], 
    		\end{equation}
set $\sigma \leftarrow  \sigma_{new}$, and go to Step~5.1.\\

\noindent
{\bf Remark }
 In Assumption~A2 we will specify the way in which we choose $z^{trial}$ in (\ref{subproblemfeas}).

\section{Plan of the proofs} \label{purpose}

  The goal of the present research is to show that, using BIRA and  under suitable assumptions, 
convergence to 
   feasible and optimal solutions takes place and worst-case 
   complexity results, that provide bounds on the evaluation 
    computer work used by the algorithm in terms
   of given small tolerances, can be proved.  These results will be
   stated in Theorems~\ref{teocomplexidade}  and \ref{teoconvergencia}. 
   
   The main assumption 
    in these theorems is that the algorithm does not stop by Restoration Failure. 
  Note that the possibility
   of stopping by Restoration Failure is unavoidable in any algorithm for constrained optimization as, 
  in some cases, feasible solutions may not exist at all. In our approach   
   optimality will be  measured by means of the Euclidean projection of the gradient of the objective 
   function onto the tangent approximation to the constraints. This is related to using 
   the Sequential Optimality Condition called L-AGP in \cite{ahm}. Such condition holds 
   at a  local minimizer of constrained optimization problems without invoking constraint qualifications. 
    Under weak constraint qualifications, the fulfillment of L-AGP implies KKT conditions \cite{amrs}.

Let us draw, now,  the general map along which the main results of the paper are proved.

    \begin{enumerate}
\item The success of the method proposed in this work is associated with the decrease of the infeasibility $\|h(x, y)\| + g(y)$,
  that should go to zero, and the decrease  of the merit function, which, ultimately, should behave as the true objective
   function onto the feasible region.
\item The iteration of the main algorithm BIRA begins calling
    Algorithm RESTA, which forces the improvement of similarity (precision) and feasibility of algebraic constraints $h(x, y)=0$. However, 
   RESTA may fail because the original problem could be infeasible. In this case BIRA stops.
\item At each iteration $k$, after success of RESTA, we update the penalty parameter $\theta$ that defines
   the merit function and we go to the Optimization Phase. At the first $N_{acce}$ attempts of the 
   Optimization Phase we try to improve optimality  without necessarily
  increasing precision in evaluations. For example, it is interesting, in practical implementations, to try  
   $y^{k+1}=y^k$ at the first iterations of the Optimization Phase. If we are not successful in the first  $N_{acce}$ attempts,
 we improve the precision taking $y^{k+1}=\yrest$, as computed by RESTA. In any case, given the accuracy level induced by $y^{k+1}$, we try to improve optimality using quadratic programming, and we test if 
  sufficient decrease of both the objective and the merit function were obtained. If this is the case, the iteration finishes.
   Otherwise, the regularization parameter that defines the quadratic programming problem is increased and quadratic minimization
    is employed again.
\end{enumerate}

The description given above induces the map of the proofs presented in this paper. Firstly, we need to prove that each iteration
  is well defined. Looking at the steps described above, for
    this purpose we need to prove that RESTA is well defined and stops in finite time. This is done in Section~5.
  Moreover we need to prove that the Optimization Phase finishes in finite time
  too. This fact will be proved in Section~6.

In Section 7 we
    prove that the infeasibility measure tends to zero. This fact
  is essential to show that, ultimately, the algorithm finds solutions of the original problem.

  Finally, in Section 8 we show that, not only the infeasibility measure but also a suitable 
  optimality measure tends to zero.

 In all the cases, convergence results are complemented with complexity results. That is, we will prove not only that crucial
  quantities produced by the algorithm tend to zero, but also that the computer work necessary to reduce those quantities to 
  a small tolerance is suitably bounded as a function of the tolerance.

\section{General Assumptions} \label{assumptions}
 
The assumptions stated in this section are supposed to hold all along this paper without 
 specific mention. These assumptions state regularity and boundedness of the functions involved
  in the definition of the problem. 

\noindent
{\bf Assumption G1 } Differentiability of $f$: The function 
  $f(x,y)$ is continuously differentiable with respect to $x$ for all 
  $x \in \Omega$  and all $y \in Y$.

\noindent
{\bf Assumption G2  } Boundedness: There exists 
  $ C_f > 0$ such that, for all $ x \in \Omega $ and for all  $y \in Y$, we have that 
 \begin{equation}\label{boundf}
 	|f(x,y)| \leq   C_f. 
 \end{equation}

\noindent
{\bf Assumption G3 } Lipschitz-continuity: There exists 
 $  L_f  \geq 0 $ such that, for all  $ x_1, x_2 \in \Omega $ and all $ y \in Y $, we have that:
 \begin{equation}\label{lipsf}
 	|f(x_1,y) - f(x_2,y)|  \leq  L_f \norm{ x_1 - x_2}
 \end{equation}
 and    
 \begin{equation}\label{lipsgradf}
 	\norm{ \nabla_x f(x_1,y) - \nabla_x f(x_2,y) }   \leq  L_f \norm{x_1 - x_2 }. 
 \end{equation}

\noindent
{\bf Assumption G4 } Upper bound for $f$: For all 
  $ x_1, x_2 \in \Omega $  and all $ y \in Y $ we have that
 \begin{equation}\label{ftaylor}
 	f(x_2,y)  \leq  f(x_1,y) + \nabla_x f(x_1)^T (x_2 - x_1)  + L_f \norm{x_2-x_1}^2.
 \end{equation}

\noindent
{\bf Assumption G5 } Differentiability of $h$: The function $h(x, y)$ is continuously differentiable with respect to $x$ for all 
  $x \in \Omega$  and all $y \in Y$.

\noindent
{\bf Assumption G6 } Boundedness of 
  $\norm{h}$  and $\norm{\nabla_x h}$: There exists $ C_h \geq 0$ such that, for all $ x \in \Omega $ and all $ y \in Y $, we have that
 \begin{equation}\label{boundu} 
 	\norm{h(x,y)} \leq C_h 
 \end{equation}
   and
 \begin{equation}\label{boundgrau}
 	\norm{ \nabla_x h(x,y) } \leq C_{h}.
 \end{equation}

\noindent
{\bf Assumption G7 } Lipschitz-continuity of $h$ and $\nabla_x h$: There exists $L_h \geq 0$ such that, for all
 $ x_1, x_2 \in \Omega $ and all  $ y \in Y $, we have that:
 \begin{equation}\label{lipsu}
 	\norm{h(x_1,y) - h(x_2,y)} \leq L_h \norm{x_1 - x_2}, 
 \end{equation}
 and  
 \begin{equation}\label{lipsgradu}
 	\norm{ \nabla_x h(x_1,y)^T - \nabla_x h(x_2,y)^T} \leq L_h \norm{x_1 - x_2}.  
 \end{equation}
 
\noindent
{\bf Assumption G8 } Upper bound of $\norm{h}$:  For all 
  $x_1, x_2 \in \Omega$ and all $y \in Y$ we have that
 \begin{equation}\label{utaylor}
 	\norm{h(x_2,y)}  \leq \norm{h(x_1,y) +  \nabla_x h(x_1,y)^T (x_2-x_1) } + L_h \norm{x_2 - x_1}^2.
 \end{equation}

\noindent
{\bf Assumption G9} Boundedness of $g_f$ and $g_h$: There exists 
  $C_g \geq 1$ such that 
 \begin{equation}\label{boundg} 
 	g_f(y) \leq C_g \mbox{ and } g_h(y) \leq C_g  
 \end{equation} 
for all $y \in Y$. 
 
\noindent
{\bf Assumption G10} Differentiability of $c(x, y)$: The function $c(x, y)$, defined by (\ref{defc}),   
is continuously differentiable with respect to $x$ for all  
 $x \in \R^n $ and $y \in Y$.

\noindent
{\bf Assumption G11} Lipschitz continuity of 
  $\nabla_x c$: There exists $ L_c \geq 0$ such that for all $ x_1, x_2 \in \Omega $ and all $ y \in Y $, we have that 

 \begin{equation}\label{lipsgradc}
 	\norm{\nabla_x c(x_1,y) - \nabla_x c(x_2,y) } \leq L_c \norm{x_1 - x_2}. 
\end{equation}

\noindent
{\bf Assumption G12} Upper bound of $c(x, y)$: For all  
 $x_1, x_2 \in \Omega$ and all  $y \in Y$ we have that
 \begin{equation} \label{ctaylor}
 	c(x_2,y) \leq c(x_1,y) + \nabla_x c(x_1,y)^T (x_2 - x_1)  + L_c \norm{x_2-x_1}^2.
 \end{equation}



\section{Theoretical Results Concerning the Restoration Phase} \label{restauracao}

   The Restoration Phase is the subject of Step~1 of BIRA. This phase begins acknowledging
  the possibility that, using some problem-dependent procedure (PDP), one may be able 
   to compute $x^k_R$ and $y^k_R$ fulfilling the conditions 
  (\ref{erre3}) and (\ref{beta3}).

    If there is no problem-dependent procedure that computes $x^k_R$ and $y^k_R$
  satisfying (\ref{erre3}) and (\ref{beta3}) we try improve feasibility executing 
    steps 2--5 of  RESTA.  However, even Algorithm RESTA may fail in that purpose, and in this case
   we declare ``Restoration Failure" and Algorithm BIRA stops. Note that every algorithm 
   for constrained optimization may fail to find feasible points, unless special conditions 
   are imposed to the problem. The main reason is that, in extreme cases, feasible points 
   could not exist at all. 

   The idea of RESTA is to show that, using quadratic programming resources, we are able to
    compute a condition similar to (\ref{beta3}). This means that only (\ref{fallar}) may fail 
    to occur in cases of probable infeasibility. 

    Assumption~\ref{restfinita} states that finding $w^{i+1}$ at Step~2 of RESTA is always inexpensive.
   The reason is that, in general, (\ref{restorationg}) merely represent increasing the precision
   in which the objective function and the constraints will be evaluated. \\

   \begin{assump}\label{restfinita}
   Step~2 of Algorithm RESTA, leading to the definition of $w^{i+1}$ satisfying   
  (\ref{restorationg}), 
    can be computed in finite time for all $k$ and $i$, without evaluations of $f$ or $h$.
\end{assump}

    At Step~5.1 of RESTA we defined $z^{trial}$ as an approximate solution of 
   problem (\ref{subproblemfeas}). Assumption~\ref{reduzc} states a simple condition that 
   such approximate solution must satisfy. According to this very mild assumption, the 
    trial point $z^{trial}$ should not be worse than $z^\ell$ in terms of functional value. 
   Note that even $z^{trial} = z^\ell$ satisfies this assumption. \\

  \begin{assump}\label{reduzc}
	For all $z^{\ell}$ and $w^{i+1}$, the point $z^{trial}$ found
 at Step~5.1 of Algorithm RESTA satisfies:

	\begin{equation} \label{eqreduzc} \nabla_x c(z^\ell,w^{i+1})^T (z^{trial} -  z^\ell) 
   + \frac{1}{2} (z^{trial} - z^\ell)^T ( B_\ell+\sigma I) (z^{trial} - z^\ell)  \leq 0.
	\end{equation} 
\end{assump}

   In Lemma~\ref{restorationcontrolc} we prove that, taking the regularization parameter 
  $\sigma$ large enough when solving (\ref{subproblemfeas}) we obtain sufficient reduction 
   of the sum of squares infeasibility at the approximate solution $z^{trial}$. In other words, 
   the loop at Steps~5.1--5.2 of RESTA necessarily finishes with the fulfillment of (\ref{desc}).\\

 \begin{lemma} \label{restorationcontrolc} 
	Suppose that Assumptions \ref{restfinita} and \ref{reduzc} hold. Define
		$\bar{\sigma} = 2 \left( L_c +\frac{M}{2}+ \alpha_R\right). $
  Then, if
  $ z^{trial} $ is computed at Step~5.1 with  $ \sigma \ge \bar{\sigma} $,  we have that 
	\begin{equation*}
		c(z^{trial},w^{i+1}) \leq c(z^\ell,w^{i+1}) - \alpha_R \|z^{trial}- z^\ell\|^2.
	\end{equation*}
     As a consequence, for all
  $k$, $i$, and $\ell$ we have that  $\sigma \leq \max\{10 \bar{\sigma}, \sigma_{max}\}$. 
\end{lemma}

  Assumption \ref{assumpprojsubprobrest} adds an additional condition that must be satisfied by  the approximate 
  solution of the subproblem (\ref{subproblemfeas}). It will be required that 
    an approximate optimality condition, expressed in terms of the projected gradient of the 
   objective function of (\ref{subproblemfeas}), should be fulfilled with tolerance proportional 
   to $\|z^{\ell+1}-z^\ell\|$.\\

 \begin{assump}\label{assumpprojsubprobrest}
  There exists $\kappa_R > 0$ such that, whenever  $ z^{\ell + 1} $ is defined at Step~5.1 of RESTA, we have that:
	\begin{equation}\label{projnula}\begin{array}{c}
\|P_\Omega \left( z^{\ell+1} - \left[ \nabla_x c(z^\ell,w^{i+1}) + B_\ell(z^{\ell+1} - z^\ell) + \sigma (z^{\ell+1} - z^\ell) \right] \right)- z^{\ell+1} \| \\\leq \kappa_R \|z^{\ell+1} - z^\ell\|.
	\end{array}
	  \end{equation} 
\end{assump}    

As a consequence of the previous assumptions, Lemma~\ref{thefo} proves that the projected gradient 
  of the linear approximation of the sum of squares at $z^\ell$, computed at the subproblem 
  solution $z^{\ell+1}$, is proportional to the norm of the difference between $z^\ell$ and 
   $z^{\ell+1}$. \\       

 \begin{lemma} \label{thefo}
   Suppose that Assumptions A1, A2, and A3 hold. Define 
		$c_{P_\Omega} =   L_c + M+\kappa_R+\max\{10\bar{\sigma}, \sigma_{max}\},$
 where $\bar{\sigma}$ was defined in Lemma \ref{restorationcontrolc}.              
Then,  
  whenever $ z^{\ell + 1} $ is defined at Step~5.2 of RESTA,  we have:
	\begin{equation} \label{graleq0}
		\norm{P_\Omega \left( z^{\ell+1} - \nabla_x c(z^{\ell+1},w^{i+1}) \right) - z^{\ell+1} } \leq  c_{P_{\Omega}}   \|z^{\ell+1}-z^\ell\|.   
	\end{equation}
\end{lemma}

   Lemma~\ref{resumoGencan} 
  establishes that, in a bounded finite number of steps, the sum of squares of 
   infeasibilities is smaller than $   r^2  c(x^k,w^{i+1})    $ or its projected gradient at $z^\ell$ is smaller than
   $  r_{feas} \|h(x^k,w^{i+1})\| $.  In other words either the squared residual or its projected gradient is smaller 
   than a multiple of the residual norm at the current iterate.   \\

 \begin{lemma} \label{resumoGencan}
   Suppose that Assumptions A1, A2, and A3 hold. Define 
 $  C_{rest}= \frac{c_{P_\Omega}^2(1-r^2)}{2\alpha_R\,  r_{feas}^2}+1,$
where $c_{P_\Omega}$ is defined in Lemma \ref{thefo}. 
 Then, at every call to Algorithm~RESTA, there exists $\ell \leq C_{rest}$ such that, defining
      	\begin{equation}\label{parafeas}
    	c_{target}  = r^2  c(x^k,w^{i+1}) \quad 
    	\text{ and }  \quad    \epsilon_c  =   r_{feas} \|h(x^k,w^{i+1})\| .
    	\end{equation} 
 we have that 
	\begin{equation}\label{restctargetz}
		c(z^{\ell},w^{i+1}) \leq c_{target}
	\end{equation}
	or
	\begin{equation}\label{restnablacz}
		\norm{ P_{\Omega}(z^{\ell} -  \nabla_x c(z^\ell,w^{i+1}))-z^{\ell} } \leq \epsilon_c.
	\end{equation}    
\end{lemma}

  The following is a technical assumption that involves $z^{\ell+1}$ obtained at Step~5 of 
   RESTA. It states that, if we add to (\ref{subproblemfeas}) the constraint that 
   $z- z^\ell$ is a multiple of $z^{\ell+1}-z^\ell$, the corresponding solution is 
   close to $z^{\ell+1}$. Clearly, this assumption holds if $z^{trial}$ is the global
  solution of (\ref{subproblemfeas}) and very plausibly holds for approximate solutions. \\

  \begin{assump}\label{assumpdistminunidimen}
   There exists $\kappa_\varphi>0$ such that,  whenever  
	$ z^{\ell + 1} $ is  the approximate solution of 
 \eqref{subproblemfeas} obtained in RESTA and 
 $z^{\ell+1}_*$ is an exact solution to the problem that has the same objective
   function and constraints as (\ref{subproblemfeas}) and, in addition 
   a constraint saying that  
  $z-z^{\ell}$ is a multiple of   $ z^{\ell + 1}-z^{\ell}$, we have:
	\begin{equation}\label{eqdistminunidimen}
		\| z^{\ell+1} -z^\ell\| \leq \kappa_{\varphi} \|z^{\ell+1}_* - z^\ell\|.
	\end{equation} 
\end{assump}

   In the following lemma we prove that the difference between consecutive internal 
   iterations in RESTA is proportional to the infeasibility at $x^k$.\\

 \begin{lemma}  \label{lemascurto}
 Suppose that Assumptions A1--A4 hold. Define
\begin{equation} \label{ceese}
 C_s  =  \kappa_{\varphi} M C_{h},
\end{equation}
where $C_h$ is defined in Assumption G6. 
Then, for all $k$, $i$ and $\ell$, 
 the iterates generated in RESTA satisfy 
	\begin{equation}\label{cotazeta0}
		\|z^{\ell+1} - z^\ell\|   \leq C_{s}  \|h(x^k,w^{i+1})\|.
	\end{equation} 
\end{lemma}

   In Lemma~\ref{Navhrest} we prove that, at every call of RESTA, the descent condition 
  on the sum of squares of infeasibilities (\ref{desc}) is tested a finite number of times.\\ 

  \begin{lemma} \label{Navhrest}
   Suppose that Assumptions A1--A4 hold. 
  Define 
		$n_\sigma =  \lfloor \log_2 (\bar{\sigma})- \log_2(\sigma_{min})\rfloor +1 $
and 
	$	N_{RESTA} = (C_{rest} n_\sigma+1) N_{prec}.$
Then, at every call to RESTA,  
  the number of tests of the condition \eqref{desc}  and the number of
  evaluations of  $h$ and $\nabla_x h $ is bounded by $N_{RESTA}$.
\end{lemma}

 In the following lemma we prove that the norm of the difference between the restored point
  $x^k_R$ and the current point $x^k$ is bounded by a multiple of $\|h(x^k, y^k_R)\|$. \\

 \begin{lemma} \label{valecondbetasobreeps} 
 Suppose that Assumptions A1--A4 hold. Define 
	$	\beta_R = \linebreak \max\{ \beta_{PDP},\beta_c+N_{RESTA} C_{s} \}$, 
where $C_s$ is defined by (\ref{ceese}). 
Then, for every iteration $k$ of BIRA, 
$(x_R^k,y_R^k)$ satisfies 
	\begin{equation} 
		\norm{ x_R^k-x^k} \leq {\beta_R}  \|h(x^k,y^k_R)\| .
	\end{equation}
\end{lemma}

  In Lemma~\ref{boudfbyinfeas} we prove that the deterioration in the objective 
   function in $x^k_R$ with respect to the objective function at $x^k$ is bounded by 
   quantity that is proportional to the infeasibilities $h$ and $g$.  For proving that result 
  we need a final assumption that states that fixing $x^k$ and restoring $y^k$ the deterioration
   in $f$ 
   is smaller than a multiple of $g(y^k)$.\\

   \begin{assump}\label{restgbeta}
 There exists $\beta >0$ such that, for all iteration $k$, ,  
	\begin{equation} \label{restorationf}
		f(x^k,y^k_R)
  \leq f(x^k,y^k) + \beta g(y^k)
  .
	\end{equation}
\end{assump}

\begin{lemma} \label{boudfbyinfeas} 
 Suppose that Assumptions A1--A5 hold. Define 
	$	\beta_f=  L_f\beta_R+ \beta.$
 Then, for every iteration $k$ of Algorithm BIRA, the point 
 $(x_R^k,y_R^k)$ computed at Step~5 of the Restoration Phase, satisfies        
	\begin{equation*}
		f(x_R^k,y^k_R) \le f(x^k,y^k)  + \beta_f[   \norm{h(x^k,y^k_R)} +g(y^{k})].
	\end{equation*}
\end{lemma}

Finally, in Assumption~\ref{assumpinexpensive} we state the sense in which the problem-dependent restoration
  procedure PDP is considered to be inexpensive. Then, in Theorem~\ref{teofinres}, the main results
  of the present section are condensed.

 \begin{assump}\label{assumpinexpensive}
	There exists $N_{PDP} $, independent of $\bar{\epsilon}_{prec}$, such that if the  problem-dependent restoration 
procedure is used at Step~1 of RESTA, it employs at most $N_{PDP} $ evaluations of  $h$ and $\nabla_x h $ and no evaluation 
of $f$ and $\nabla_x f $.
\end{assump}

  Successful restoration procedures in IR methods usually satisfy stability conditions that say that the distance between
   restored points and current iterates is bounded by a constant times the infeasibility measure. Alternatively, 
    it is generally proved that the objective function at the restored point is smaller than the objective function at 
   the current iterate plus a constant times the infeasibility. The stability conditions obviously hold when $y^k_R$ and 
   $x^k_R$ are computed by the problem-dependent procedure PDP, as stated in (\ref{beta3}). In Theorem~\ref{teofinres} we prove
    that similar results hold in the case that restoration is achieved by means of Steps 2--5 of RESTA.\\

      \begin{theorem}  \label{teofinres} 
   Suppose that the General Assumptions and Assumptions A1--A6 hold. Then, there exist $N_R$ and $\beta_f$, independent
 of $\bar{\epsilon}_{prec}$, such that, for every iteration $k$ of Algorithm BIRA, the point  $(x_R^k,y_R^k)$ is computed employing at most $N_R$ 
   evaluations of  $h$ and $\nabla_x h $, no evaluation of $f$ and $\nabla_x f $, satisfying      
   	\begin{equation} \label{cotaxyfeas}  
   	\norm{ x_R^k-x^k} \leq {\beta_R}  \|h(x^k,y^k_R)\| .
   \end{equation}
and 
   \begin{equation}\label{cotafxyrestoration}
   	f(x_R^k,y^k_R) \le f(x^k,y^k)  + \beta_f[   \norm{h(x^k,y^k_R)} +g(y^{k})].
   \end{equation}
    \end{theorem}  

\begin{proof}
	Conditions \eqref{cotaxyfeas} and \eqref{cotafxyrestoration} follow 
 from  Lemmas \ref{valecondbetasobreeps} and  \ref{boudfbyinfeas}, respectively. Observe that 
 no evaluation of $f$ and $\nabla_x f $ is made when calling {Algorithm}~RESTA. So, defining $N_R=N_{RESTA}+N_{PDP}$, 
by Lemma   \ref{Navhrest} and Assumption~\ref{assumpinexpensive}  we have the desired result.  
\end{proof}

\section{BIRA is well defined} \label{birawell}

All along this section we will  assume, without specific mention,
    that the General Assumptions G1--G12 and the
  Restoration Assumptions A1--A6 are fulfilled. Assumption~\ref{bajasub} will be added when needed
   to prove specific results and its fulfillment will be mentioned whenever necessary. 

  As the title of this section indicates, the objective will be that Algorithm BIRA is well defined, 
   that is, that for any iteration of BIRA, either the algorithm stops or it is possible to 
  compute the next iterate.\\

  We begin showing that the penalty parameter is well defined and satisfies the inequality
 (\ref{penaltetamas}), that states that, from the point of view of the merit function, 
   the restored point $x^k_R$ is better than the current iterate $x^k$. \\

      \begin{lemma} \label{tetabemdef} 
At every iteration~k of BIRA, the  penalty parameter $\theta_{k+1}$ is well defined, 
 $0<\theta_{k+1} \leq \theta_k$, and 
        \begin{equation} \label{penaltetamas}
        \begin{array}{c}
        \Phi(x_R^k, y_R^k,\theta_{k+1})-\Phi(x^k,y^k_R, \theta_{k+1})\\ 
\le \frac{1-r}{2}  \left[ \|h(x^k_R,y_R^k)\| - \|h(x^k,y^k_R)\| + g(y_R^k) -g(y^k)\right].
        \end{array}
        \end{equation}     
    \end{lemma}

In Lemma~\ref{lemabarteta} we prove that the penalty parameters are bounded away from zero.\\ 
 
      \begin{lemma} \label{lemabarteta} 
 Define $
         \bar{\theta}= \min\left\{\theta_0,\left[\frac{2}{1+r} \left(\frac{L_f \beta_R }{1-r} + 1\right) \right]^{-1}\right\}.
 $
 Then, for every iteration $k$ in Algorithm BIRA we have that 
        \begin{equation} \label{barteta}
            \theta_k  \geq \bar{\theta} > 0. 
        \end{equation}
    \end{lemma}

  The following assumption establishes the conditions that must be satisfied by an approximate 
   solution of \eqref{subproblemintermediario}. \\

     \begin{assump} \label{bajasub}
   There exists $\kappa_T>0$ such that, at every iteration $k$ of Algorithm BIRA, 
   the approximate solution of the quadratic programming problem \eqref{subproblemintermediario} satisfies  
            \begin{equation} \label{dessubprob}
                \nabla_x f(x_R^k,y^{k+1})^T (x - x_R^k)  +  \frac{1}{2} (x - x_R^k)^T H_k (x - x_R^k) + \mu \|x - x_R^k\|^2 \leq 0
            \end{equation}    
        and 
        \begin{equation} \label{errotan}
        		\|\nabla_x h(x_R^k,y^{k+1})^T (x-x_R^k) \| \leq \kappa_T \|x-x_R^k\|^2. 
       \end{equation}   
    \end{assump}

  In the following lemma we prove that, when in the Optimization Phase, for a sufficiently large regularization parameter $\mu$, the 
  descent conditions for the objective function and the merit function are satisfied. As a consequence,
  in the subsequent corollary we establish the maximal number of iterations that could be needed
  to fulfill those conditions.\\

     \begin{lemma} \label{lemmaNreg0} 
    Suppose that Assumption \ref{bajasub} holds. Define
    $C_\mu = M +  \tilde{\alpha}+ L_f,$
 where
     \begin{equation} \label{defalfatil}
       \tilde{\alpha}= \max\left\{ \alpha, \frac{ 1 - \bar{\theta} }{\bar{\theta}} (\kappa_T+L_h) \right\}.
    \end{equation}          
  Then, if  $ \mu \ge C_\mu$, $y^{k+1}=y^k_R$, and $x$ is the solution of \eqref{subproblemintermediario}, the conditions \eqref{desfinte}  and 
 \eqref{desmeri} are fulfilled. 
    \end{lemma}

     \begin{cor}\label{lemmaNreg}
            Suppose that Assumption~\ref{bajasub} holds. Define 
$       N_{reg} = \linebreak\max\{ \lfloor \log_2(C_\mu)-\log_2(\mu_{min})\rfloor,   N_{acce}\}+1$
   and
$\bar{\mu} =  \max\{10 C_{\mu},10^{N_{acce}} \mu_{max}\}. $
  Then, after at most $N_{reg}$ sub-iterations at Step~3 of BIRA, the conditions  \eqref{desfinte} and
 \eqref{desmeri} are fulfilled . Moreover, 
 $\mu_k \leq \bar{\mu}$ for all $k$.
    \end{cor}

 \section{Convergence to feasibility}\label{conviastep}

    In this section we will prove that, when executing BIRA, the infeasibility measure 
   tends to zero. Moreover, we will prove a crucial theorem which shows
       that the norm of the difference between $x^{k+1}$ and $x^k_R$   
     tends to zero.

 For all the proofs of this section we
    will assume, without specific mention,  that all the General Assumptions,  
 the Assumptions A1--A7, and the following Assumption~\ref{restgbeta2} take place. 
   Assumption~\ref{restgbeta2} states that bounded deterioration of objective function and 
  also $h$-feasibility occurs in a restricted way, depending on a possibly small parameter
   that depends on $\bar{\theta}$. This means that, in the worst case, bounded deterioration 
   with respect to precision does not occur at all. Note that, however, the new bounded 
   deterioration condition needs to hold only for $k$ large enough. \\

     \begin{assump}\label{restgbeta2}
     Let   $\bar{\theta}$ be as defined in Lemma   \ref{lemabarteta}. Then, there exist $k_R$, and $\gamma \in (0,1)$ such that, for $\bar{\beta} =\frac{\bar{\theta}(1-\gamma)(1-r)^2}{2}$ and all $k \ge k_R$, 
    	\begin{equation} \label{restorationf2}
    	f(x^k,y^{k+1})
    		 \leq f(x^k,y^k) + \bar{\beta} g(y^k) \quad \text{ and } \quad \norm{ h(x^k,y^k_R)  } \leq \norm{ h(x^k,y^{k})}  + \bar{\beta} g(y^k).
    	\end{equation}
    \end{assump}

   Theorem \ref{teohedvaoa0} states the summability of all infeasibilities.\\
    
     \begin{theorem}  \label{teohedvaoa0} 
      Define 
     \begin{equation}\label{defcfeas}
        C_{feas} = \frac{2}{\gamma (1-r)^2} \left[  k_R (2C_f+C_h)+ C_\rho+ C_h +C_g  \right].
    \end{equation}          
      Then, 
        \begin{equation} \label{limsomah}
            \sum_{j=0}^{k}[ \|h(x^j,y^j_R)\| + g(y^j) ] \leq C_{feas}.
        \end{equation}
    \end{theorem}

    \begin{proof}
      Let us define 
    \begin{equation}\label{defrho}
        \rho_j = \frac{1-\theta_j}{\theta_j}=\frac{1}{\theta_j}-1,  \text{ for all  $j \leq k$.}
    \end{equation}
     By Lemma  \ref{lemabarteta}, we know that $\theta_j \in (0, 1)$,   
      $\{\theta_j\}$ is non-increasing and bounded below by $\bar{\theta}$. 
     Then, the sequence 
 $\{\rho_j\}$ is positive, non-decreasing and bounded above by  $\bar{\rho}= \dfrac{1}{\bar{\theta}} -1$. So, since $\rho_0 > 0$,  
     \begin{equation}\label{suma}
        \sum_{j=0}^{k-1} (\rho_{j+1}-\rho_{j} ) = \rho_{k}- \rho_{0} <  \rho_k = \frac{1}{\theta_k} -1 \le  \frac{1}{ \bar{\theta} } -1 < \frac{1}{ \bar{\theta} } < \infty. \;   
    \end{equation}
   By \eqref{boundu}, we have that
   $\norm{h(x^j,y^{j+1}) } \leq C_h$ for all  $j$.  Since $\rho_{j+1}-\rho_{j} \geq 0$, taking $C_{\rho} \equiv \frac{C_h}{\bar{\theta}} $, thanks to \eqref{suma}, we have that   
    \begin{equation} \label{suma2}
        \sum_{j = 0}^{k-1} (\rho_{j+1} - \rho_j) \|h(x^j,y^{j+1})\| \leq\sum_{j = 0}^{k-1} (\rho_{j+1} - \rho_j)C_h  \leq  \frac{C_h}{\bar{\theta}} = C_{\rho} < \infty.
    \end{equation}

    We have that
     \begin{equation} \label{desde}
     	\begin{array}{c}
     	\Phi(x^{j+1}, y^{j+1},  \theta_{j+1}) - \Phi(x^j, y^{j+1}, \theta_{j+1})\\ \leq    \frac{1-r}{2}\left[ \|h(x_R^j,y_R^j)\| - \|h(x^j,y^j_R)\| + g(y_R^j)- g(y^j) \right]\\ 
     	\leq -   \frac{(1-r)^2}{2}\left[  \|h(x^j,y^j_R)\| + g(y^j)\right], 
     	\end{array}             
    \end{equation}
      where the second inequality comes from $\|h(x_R^j,y_R^j)\|  \leq r  
       \|h(x^j,y^j_R)\| $ and   $g(y_R^j) \le r g(y^j)$. 
    
    By the definition of $\Phi$, dividing \eqref{desde} by  
      $ \theta_{j + 1} $, we have that, for all  $j \leq k-1$,
   \begin{equation*}\begin{array}{c}
          f(x^{j+1},y^{j+1})+ \frac{1 - \theta_{j+1}}{\theta_{j+1}} \left[\|h(x^{j+1},y^{j+1})\| + g(y^{j+1})\right]  -  f(x^j,y^{j+1})   \\  -\frac{1 - \theta_{j+1}}{\theta_{j+1} } [ \|h(x^j,y^{j+1})\| +g(y^{j+1})] \leq -  \frac{(1-r)^2}{2 \theta_{j+1}}  \left[  \|h(x^j,y^j_R)\| + g(y^j) \right].\end{array} 
    \end{equation*} 
      By the definition of 
      $\rho_j$ in \eqref{defrho}, using  that $\theta_j \in (0,1)$, we deduce that
    \begin{equation} \label{e99}
    	\begin{array}{c}
 	\frac{(1-r)^2}{2}\left[  \|h(x^j,y^j_R)\| + g(y^j)\right] \leq 	\frac{(1-r)^2}{2 \theta_{j+1}} \left[  \|h(x^j,y^j_R)\| + g(y^j)\right] \\ 
 	\leq  f(x^j,y^{j+1})-f(x^{j+1},y^{j+1})+\rho_{j+1}\|h(x^j,y^{j+1})\|-\rho_{j+1}\|h(x^{j+1},y^{j+1})\|.
    	\end{array}
    \end{equation}
     Adding and subtracting 
      $\rho_j \|h(x^j,y^{j+1})\|$ on the right-hand side of (\ref{e99}), and arranging terms, we have: 
\begin{equation*}\begin{split}
        \frac{(1-r)^2}{2 } [\|h(x^j,y^j_R)\|  +g(y^j)]  & \leq  f(x^j,y^{j+1}) -  f(x^{j+1},y^{j+1})\\
        &\quad +  (\rho_{j+1} - \rho_{j}) \|h(x^j,y^{j+1})\|  \\
       &\quad   +  \rho_{j} \|h(x^j,y^{j+1})\|  -  \rho_{j+1} \|h(x^{j+1},y^{j+1})\|  . \end{split}
    \end{equation*} 

  Then,     
    \begin{equation} \label{limitainv}\begin{split}
       \frac{(1-r)^2}{2 }{ \sum_{j=0}^{k-1}} [\|h(x^j,y^j_R)\| +g(y^j)] & \leq      f(x^0,y^1) - f(x^k,y^k)\\
         &\quad +  \sum_{j=1}^{k_R-1} [f(x^j,y^{j+1}) -  f(x^{j},y^{j})]\\
         &\quad  +  \sum_{j=k_R}^{k-1} [f(x^j,y^{j+1}) -  f(x^{j},y^{j})]\\
         & \quad   +{ \sum_{j=0}^{k-1}} (\rho_{j+1} - \rho_{j}) \|h(x^j,y^{j+1})\|\\
  & \quad +  \rho_0 \|h(x^0,y^1)\|-   \rho_{k} \|h(x^k,y^k)\| \\
         & \quad + \sum_{j=1}^{k_R-1} \rho_j [\|h(x^j,y^{j+1})\| - \|h(x^{j},y^{j})\|]\\& \quad + \sum_{j=k_R}^{k-1} \rho_j [\|h(x^j,y^{j+1})\| - \|h(x^{j},y^{j})\|].
    \end{split}
    \end{equation} 
    By \eqref{suma2}, ${\sum_{j=0}^{k-1}} (\rho_{j+1} - \rho_j) \|h(x^j,y^{j+1})\| \leq C_{\rho} $. 
  Moreover, since  $\rho_j \leq \bar{\rho}$, by Assumption~\ref{restgbeta2} and disregarding the certainly non-positive terms,  \eqref{limitainv} implies that
    \begin{equation*}
    	\begin{split}
         \frac{(1-r)^2}{2 } \sum_{j=0}^{k-1} [\|h(x^j,y^j_R)\|  +g(y^j)] &   \leq  |f(x^0,y^1)|+ |f(x^k,y^k)| \\ &\quad + \sum_{j=1}^{k_R-1} |f(x^j,y^{j+1})- f(x^{j},y^{j})|\\
         &\quad + \sum_{j=k_R}^{k-1} \bar{\beta} g(y^j) +C_\rho +  \rho_{0} \|h(x^0,y^1)\|   \\& \quad +  \sum_{j=1}^{k_R-1} \bar{\rho} \|h(x^j,y^{j+1})\|+  \sum_{j=k_R}^{k-1} \bar{\rho} \bar{\beta} g(y^j). 
    \end{split}
    \end{equation*} 
    
    By \eqref{boundf}, \eqref{boundu}, and \eqref{boundg} we have that $ f $, $\norm{h} $,  and $g$  are bounded above
  by $ C_f $ , $ C_{h} $, and $C_g$ respectively. Then, as $\bar{\rho}+1=\frac{1}{\bar{\theta}}$, we obtain that
     \begin{equation*}\begin{split}
    \frac{(1-r)^2}{2 } \sum_{j=0}^{k-1}  [\|h(x^j,y^j_R)\|+g(y^j)]   & \leq    k_R (2C_f+C_h)+ C_\rho    +\frac{\bar{\beta}}{\bar{\theta}} \sum_{j=k_R}^{k-1}  g(y^j).
 \end{split}
    \end{equation*} 

  Therefore, using that 
 $0 \leq g(y^j) \leq g(y^j)+\|h(x^j,y^j_R)\|$ and $\frac{\bar{\beta}}{\bar{\theta}} = \frac{(1-\gamma)(1-r)^2}{2}$ we obtain:  
 \begin{equation*}\begin{split}
		\frac{(1-r)^2}{2 } \sum_{j=0}^{k}  [\|h(x^j,y^j_R)\|+g(y^j)]   & \leq  k_R (2C_f+C_h)+ C_\rho\\
  & \quad + \frac{(1-\gamma)(1-r)^2}{2} \sum_{j=k_R}^{k-1}  g(y^j)\\& \quad +\|h(x^k,y^k_R)\|+g(y^k)\\
		&\leq    k_R (2C_f+C_h)+ C_\rho\\
  & \quad + \frac{(1-\gamma)(1-r)^2}{2} \sum_{j=0}^{k} \|h(x^j,y^j_R)\|\\
		&\quad +\sum_{j=0}^{k} g(y^j) +\|h(x^k,y^k_R)\|+g(y^k).
	\end{split}
\end{equation*} 
   Thus, 
   \begin{equation*}\begin{split}
   		\frac{\gamma (1-r)^2}{2 } \sum_{j=0}^{k}  [\|h(x^j,y^j_R)\|+g(y^j)]    \leq   k_R (2C_f+C_h)+ C_\rho
   		+ C_h+C_g.
   	\end{split}
   \end{equation*}    
    So, by (\ref{defcfeas}), we obtain  
     \eqref{limsomah}, as desired.\\  
    \end{proof}

 The result stated in
  Theorem~\ref{teohedvaoa1} will be used in the proof of Lemma~\ref{lemasomagradprojproporcionald} which, in turn, 
   is essential for the proof of the main theorems in Section~\ref{compcon}.\\

         \begin{theorem} \label{teohedvaoa1} 
  Define 
      $   C_d = \frac{1}{\alpha}[(\beta_f+\beta) C_{feas}  + 2C_f  ].  $
   Then, 
        \begin{equation}\label{somad}
            \sum_{j=0}^k \norm{x^{j+1} - x_R^j}^2 \leq C_d .
        \end{equation}
    \end{theorem}
    
    \begin{proof}
    	   By  \eqref{desfinte} we have that  
    \begin{equation}\label{reff}
    	\begin{array}{ll}
    	  \alpha  \norm{x^{j+1} - x_R^j}^2 &\leq f(x^j_R,y_R^j) - f(x^{j+1},y^{j+1})\\
    	  &\leq f(x^j_R,y_R^j)-f(x^j,y_R^j)+f(x^j,y_R^j)-f(x^j,y^j)\\ & \quad +f(x^j,y^j) - f(x^{j+1},y^{j+1}).
    	\end{array}  \end{equation}

  For all   
      $j \leq k-1$, by  \eqref{cotafxyrestoration}, we have that
   $ f(x_R^j,y_R^j) - f(x^j,y^j_R) \le \beta_f\left[ \|h(x^j,y^j_R)\|+ g(y^j) \right]$. On the other hand,
 \eqref{restorationf} implies that $f(x^j,y^j_R)-f(x^j,y^j)\leq \beta g(y^j)$. So, 
    \begin{equation}\label{reff2}
    		\alpha  \norm{x^{j+1} - x_R^j}^2 \leq \beta_f\left[ \|h(x^j,y^j_R)\|+ g(y^j) \right]+\beta g(y^j)+f(x^j,y^j) - f(x^{j+1},y^{j+1}). 
    	\end{equation}

   Using that  $\|h(x^j,y^j_R)\| \geq 0$ and adding terms from $j$ to $k-1$, we obtain: 
    \begin{equation*}
    \begin{array}{c}
\alpha \sum_{j=0}^{k-1} \norm{ x^{j+1} -  x_R^j }^2\\  \leq (\beta_f+\beta)\sum_{j=0}^{k-1}   \left[ \|h(x^j,y^j_R)\|+ g(y^j) \right]  +\sum_{j=0}^{k-1} [f(x^j,y^j) - f(x^{j+1},y^{j+1})].
    \end{array}
    \end{equation*}
  Therefore, by \eqref{limsomah},  
   \begin{equation}\label{calculofkappa}
   	\alpha \sum_{j=0}^{k-1} \norm{ x^{j+1} -  x_R^j }^2  \leq (\beta_f+\beta) C_{feas}+f(x^0,y^0) - f(x^k,y^k). \end{equation}
    Finally, by  \eqref{boundf} and \eqref{calculofkappa}, the desired 
  result is obtained. \end{proof}

\section{Complexity and  Convergence}\label{compcon}

In this section we suppose, without specific mention, that the General Assumptions,  
 Assumptions A1--A8,  and the following Assumption~\ref{assumpprecsubprob} hold.
  Assumption~\ref{assumpprecsubprob} merely states the approximate optimality conditions that 
    the approximate solutions of  (\ref{subproblemintermediario})  
   must fulfill. \\

\begin{assump}\label{assumpprecsubprob} 
	There exists $\kappa > 0$ such that, for every iteration   $k$ at Algorithm BIRA, 
  the approximate solutions of \eqref{subproblemintermediario}  satisfy 
	\begin{equation} \label{projectionintermed}
 \begin{array}{c}
   \norm{ P_{D^{k+1}} (x^{k+1} - \nabla_x f(x_R^k,y^{k+1}) - H_k  (x^{k+1}-x_R^k) - 2 \mu_k  (x^{k+1}-x_R^k) ) - x^{k+1} } \\  \leq \kappa \norm{  x^{k+1}-x_R^k},
 \end{array}
	\end{equation}
where  $D^{k+1}$ is defined by
\begin{equation} \label{djmas1}
	D^{k+1} = \{x \in \Omega \;|\; \nabla_x h(x_R^k,y^{k+1})^T(x-x_R^k) =0 \}.
\end{equation}
\end{assump}

In Lemma~\ref{lemagradprojproporcionald} we prove that the projected gradient of the objective 
   function onto the tangent set to the constraints tends to zero proportionally to the 
  norm of the difference between $x^{k+1}$ and the restored point $x^k_R$. \\

\begin{lemma} \label{lemagradprojproporcionald}
  Define 
\begin{equation} \label{cepe}
 C_p=  M+\kappa + 2  \bar{\mu}+2,
\end{equation} 
  where $\bar{\mu}$ is defined in Corollary~\ref{lemmaNreg}. 
 Then, 
  
	\begin{equation}
		\| P_{D^{k+1}} (x^{k}_R - \nabla_x f(x^{k}_R,y^{k+1}) ) - x^{k}_R \|
   \leq  C_p \| x^{k+1}-x_R^k  \|. \label{gradprojcontoladointermed}
	\end{equation} 
\end{lemma}

  Lemma~\ref{lemasomagradprojproporcionald} establishes the summability of squared norms of 
  the projected gradients of the objective function computed as the restored iterates.\\
 
    \begin{lemma} \label{lemasomagradprojproporcionald} Define
     $  C_{proj} = C_p^2 C_d.$
   Then, for every iteration $k$ of BIRA, we have that
 \begin{equation} \label{sumap2}
 	\sum_{j = 0}^k  \| P_{D^{j+1}} (x^{j}_R - \nabla_x f(x^{j}_R,y^{j+1}) ) - x^{j}_R \|^2 \leq   C_{proj}.
 \end{equation}
    \end{lemma}

  Lemma~\ref{lemaNfeas} is a complexity result establishing that the number of iterations at 
   which infeasibility takes place with respect to given precisions is, in the worst case, 
   proportional to the multiplicative inverse of the precisions required. From a practical point of view, to be consistent with the Restoration Failure criterion, the accuracy with respect to $g$ should be
 less demanding than the one used  in RESTA. However this is not a mathematical requirement and is not used in the following lemma. \\
 
     \begin{lemma} \label{lemaNfeas} 
    Let $\epsilon_{feas}>0$ and $\epsilon_{prec}>0$ be given. Let $N_{hinfeas}$
 be the number of iterations of BIRA at which $ \|h(x^k_R,y^k_R)\|>  \epsilon_{feas}$,
    $N_{ginfeas}$ the number of iterations of BIRA at which 
  $g(y^k)>  \epsilon_{prec}$, and $N_{infeas}$ the number of iterations of BIRA such
  that  $ \|h(x^k_R,y^k_R)\|>  \epsilon_{feas}$ or $ g(y^k_R)>  \epsilon_{prec}$.
 Then,  
        \begin{equation}\label{Nhmaiorepsf}
        \begin{split}
           &N_{hinfeas} \leq   \floor{\frac{rC_{feas} }{\epsilon_{feas}}}, \,\,\,\, N_{ginfeas} \leq     \floor{\frac{C_{feas} }{\epsilon_{prec}}}\\
           &\text{ and } \,\,\, N_{infeas} \leq   \floor{ \max\left\{\frac{rC_{feas} }{\epsilon_{feas}}, \frac{rC_{feas} }{\epsilon_{prec}}\right\}}.
           \end{split}
        \end{equation}  
     \end{lemma}

  Lemma~\ref{lemaNopt} is a complexity result that states that the number of iterations at which 
  the projected gradient of the objective function at the restored points is bigger than 
   a given tolerance $\epsilon_{opt}$ is proportional, in the worst case, to $\epsilon_{opt}^{-2}$. 

     \begin{lemma} \label{lemaNopt}  
   Suppose that $\epsilon_{opt} > 0$.
 Let $N_{opt}$ be the number of iterations such that
   $ \| P_{D^{j+1}} (x^{j}_R - \nabla_x f(x^{j}_R,y^{k+1}) ) - x^{j}_R \| > \epsilon_{opt}. $ Then, 
    \begin{equation} \label{nopt}
       N_{opt}\leq   \floor{\frac{C_{proj}}{\epsilon_{opt}^{2}}}. 
    \end{equation}
  \end{lemma}

     \begin{theorem}  \label{teocomplexidade} 
   Suppose that the General Assumptions and Assumptions A1--A9 hold. Given  
     $\epsilon_{prec} > 0$, $\epsilon_{feas} > 0$ , and $\epsilon_{opt}>0$, then:  
    \begin{itemize}
        \item If RESTA does not stop by Restoration Failure and $N_{max}$ is the maximum number of iterations $j$ of BIRA
  such that  $g(y^{j}_R) > \epsilon_{prec} $, or
  $g(y^{j+1}) > \epsilon_{prec} $, or $\|h(x^j_R,y^j_R)\|> \epsilon_{feas} $ or $\| P_{D^{j+1}} (x^{j}_R - \nabla_x f(x^j_R, \linebreak y^{j+1}) ) - x^{j}_R \| > \epsilon_{opt}$, then 
        \begin{equation}\label{Nmax}
          N_{max} \leq  \floor{\max\left\{\frac{rC_{feas} }{\epsilon_{feas}}, \frac{rC_{feas} }{\epsilon_{prec}}\right\}}+ \floor{\frac{C_{feas} }{\epsilon_{prec}}}+   \floor{\frac{C_{proj}}{ \epsilon_{opt}^{2}}}.
        \end{equation}
       
        \item The total number of evaluations of 
   $h$,  $\nabla_x h$, $f$, and $\nabla_x f$ in BIRA until declaring Restoration Failure or 
  finding $x^j_R$ such that   
        \begin{equation}\label{parada}
        	\begin{array}{c}
        		\|h(x^j_R,y^j_R)\|\leq \epsilon_{feas}, \,\,\, g(y^{j}_R) \leq \epsilon_{prec}, \,\,\, g(y^{j+1}) \leq \epsilon_{prec} \,\,\, \text{and }\\
        		\| P_{D^{j+1}} (x^{j}_R - \nabla_x f(x^j_R,y^{j+1}) ) - x^{j}_R \|\leq \epsilon_{opt}
        	\end{array}        
        \end{equation} 
    is bounded by   $N_{av}$, where 
\begin{equation}
N_{av} =O(\min\{\epsilon_{prec},\epsilon_{feas}\}^{-1}+\epsilon_{prec}^{-1}+\epsilon_{opt}^{-2}).
\end{equation}
   \end{itemize}    
    \end{theorem}  

 \begin{proof} 
  Assume firstly that BIRA does not stop with Restoration Failure. By  
  Lemma  \ref{lemaNfeas}, the inequalities
  $  \|h(x^j_R,y^j_R)\| > \epsilon_{feas}$ or $g(y^j_R)  > \epsilon_{prec}$ may occur at most during
   $\floor{\max\left\{\frac{rC_{feas} }{\epsilon_{feas}}, \frac{rC_{feas} }{\epsilon_{prec}}\right\}}$ iterations. 
 Therefore, after $\floor{\max\left\{\frac{rC_{feas} }{\epsilon_{feas}}, \frac{rC_{feas} }{\epsilon_{prec}}\right\}}+ \floor{\frac{C_{feas} }{\epsilon_{prec}}}+   \floor{\frac{C_{proj}}{ \epsilon_{opt}^{2}}} +1$  iterations, we know that at least at $\floor{\frac{C_{feas} }{\epsilon_{prec}}}+\floor{\frac{C_{proj}}{ \epsilon_{opt}^{2}}} +1$ 
  of these iterations, 
the inequalities  $g(y^j_R) \leq \epsilon_{prec}$ and $h(x^j_R,y_R^j) \leq \epsilon_{feas}$ 
 are fulfilled. 
   
   By Lemma \ref{lemaNopt}, the inequality
  $ \| P_{D^{j+1}} (x^{j}_R - \nabla_x f(x^j_R,y^{j+1}) ) - x^{j}_R\| >  \epsilon_{opt}$ 
 may occur at most in 
 $   \floor{C_{proj}/ \epsilon_{opt}^{2}} $ iterations.
 Thus, at least in 
 $ \floor{\frac{C_{feas} }{\epsilon_{prec}}}+1$ over
 $ \floor{\frac{C_{feas} }{\epsilon_{prec}}}+\floor{\frac{C_{proj}}{ \epsilon_{opt}^{2}}} +1$ 
 iterations we should have that
 $\|h(x^j_R,y^j_R)\|\leq \epsilon_{feas} $,
 $g(y^{j}_R) \leq \epsilon_{prec} $, and 
 $\| P_{D^{j+1}} (x^{j}_R - \nabla_x f(x^j_R,y^{j+1}) ) - x^{j}_R \|\leq \epsilon_{opt}$.
   
Analogously, by Lemma \ref{lemaNfeas}, the number of iterations at which
 $g(y^{j+1}) \leq \epsilon_{prec}$ is bounded by
  $\floor{\frac{C_{feas} }{\epsilon_{prec}}}$, then at least in one over the
  $ \floor{\frac{C_{feas} }{\epsilon_{prec}}}+1$ iterations \eqref{parada} takes place.
 So, \eqref{Nmax} is proved. 

Thinking   \eqref{parada} as a ``stopping criterion"
  for BIRA, the total number of iterations would be 
 at most  $N_{max}+1$, since we would have stopped by Restoration Failure  or the conditions \eqref{parada} 
would be satisfied. Let us now analyze the number of functions evaluations at  each iteration.

For every iteration of BIRA, by Lemma
\ref{Navhrest}, the Restoration Phase finishes after at most  
$N_R$ evaluations of $h$ and $\nabla_x h$. Moreover, $f$ and $\nabla_x f$ are not evaluated in the Restoration Phase.
At Step~2 of BIRA, we have two evaluations of $f$ and additional evaluations of $h$ are not necessary, 
since $h(x^k,y^k_R)$ and $h(x^k_R,y^k_R)$ have been already computed in RESTA or to check \eqref{erre3}. Furthermore, no derivatives are used  at Step~2.

Now, let us see  what happens at the Optimization Phase. Firstly, note that, for building subproblem \eqref{subproblemintermediario}, we only use one evaluation of \linebreak
$\nabla_x h(x^k_R,y^{k+1})$ and $\nabla_x f(x^k_R,y^{k+1})$  in the first $N_{acce}$ attempts of the Optimization Phase (when $y^{k+1}$ does not need to be $y^k_R$) and an extra computation of them for the remaining ones.

In the test of (\ref{desfinte}) there   is no need to calculate $f(x_R^k,y^{k}_R)$, which has already been evaluated in the Restoration Phase.  However, we need to compute
   $f(x, y^{k+1})$, which can be used for every verification of 
   \eqref{desmeri}. In every loop of the Optimization Phase it is necessary an evaluation of $h(x, y^{k+1})$ too. By Corollary~\ref{lemmaNreg}, the Optimization Phase finishes   after at most $N_{reg} $ calls to Step~3. Then,  $f(x, y^{k+1})$ and $h(x, y^{k+1})$ are evaluated at most   $ N_{reg}$ times at each call of the Optimization Phase. The values of $h(x^k_R, y^{k}_R)$ and $h(x^k, y^{k}_R)$ involved in \eqref{desmeri} have already been  computed at Step~2. Also to check  \eqref{desmeri}, we need  to compute   $f(x^k, y^{k+1})$ and $h(x^k, y^{k+1})$  once for the  $N_{acce}$ first attempts, and no additional evaluation is needed in the other ones, since  $f(x^k, y^{k}_R)$ and $h(x^k, y^{k}_R)$  have been already computed at Step~2.

 Therefore, the number of evaluations of  $h$ and $\nabla_x h$, at each iteration of 
  BIRA is, respectively, $N_R+N_{reg}+1$ and  $N_R+2$. 
  Moreover, $ f $ is computed at each iteration of BIRA at most 
  $ N_{reg} +3$ times and, at most, two evaluations of 
  $\nabla_x f$ are necessary. Since  $N_{reg}$ and $N_R$ do not depend on $\epsilon_{prec}$, $\epsilon_{feas}$, and $\epsilon_{opt}$, the total number of iterations and evaluations of $h$, $\nabla_x h$,  $ f $,  and
  $\nabla_x f$ before declaring Restoration Failure or  obtaining  \eqref{parada} is $O(\min\{\epsilon_{prec},\epsilon_{feas}\}^{-1}+\epsilon_{opt}^{-2})$.  
\end{proof}
 
 Our last theorem concerns the asymptotic convergence of  BIRA. For this, it is natural to consider that
 the algorithm generates an infinite sequence, not meeting a stopping criterion. Therefore, it is reasonable to think that $\bar{\epsilon}_{prec}$ is null. 
 
     \begin{theorem} \label{teoconvergencia}
    Suppose that the General Assumptions, Assumptions A1--A9  hold, 
     and BIRA does not stop by Restoration Failure. Then, 
    \begin{equation} \label{converge}
    	\begin{split}
    		 \lim_{k \to \infty}   g(y^k)  = 0,
    		\;\;\;\;   \lim_{k \to \infty}   g(y^k_R)  = 0,
    		\;\;\;\; \lim_{k \to \infty}  \|h(x^k_R,y^k_R)\|  = 0,\\ \text{ and  }  \;\;\;\; \lim_{k \to \infty}   \| P_{D^{k+1}} (x^k_R - \nabla_x f(x^{k}_R,y^{k+1}) ) - x^k_R \| = 0.
    	\end{split}
    \end{equation}
   
    \end{theorem}
 
    \begin{proof} Assume, by contradiction that BIRA computes infinitely many iterations
  and at least one of the sequences $\{\|h(x^k_R,y^k_R)\|\}$, $\{g(y^k)\}$ or
 $\{\| P_{D^{k}} (x^{k}_R - \nabla_x f^{k}(x^{k}_R,y^{k+1} )) - x^{k}_R \| \}$ does not
  converge to zero. To fix ideas, suppose that $\{g(y^k)\}$ does not 
  converge to zero. Then, there exists  $\varepsilon > 0$ and infinitely many indices
  $ K $ such that  $g(y^j) > \varepsilon$ for all $j \in K$. 
  Therefore,  $g(y^k) > \epsilon_{feas} $  occurs 
  infinitely many times if we define  $\epsilon_{prec} = \varepsilon$. By Theorem~\ref{teocomplexidade},
 this is impossible and so $\lim_{k \to \infty}   g(y^k)  = 0$. Since $0 \leq g(y^k_R) \leq g(y^k)$, we also have that $\lim_{k \to \infty}   g(y^k_R)  = 0$.
 
 The convergence to zero of the sequences 
    $\{\|h(x^k_R,y^k_R)\|\}$ and  $\{\| P_{D^{k}} (x^{k}_R - \nabla_x f^{k}(x^{k}_R,y^{k+1} )) - x^{k}_R\| \}$ 
 is proved in an entirely analogous way using 
 $\epsilon_{feas}=\epsilon$ or $\epsilon_{opt}=\epsilon$, respectively.
   
   \end{proof}

\section{Conclusions}  \label{conclusions}

  Many practical problems require the minimization of functions that are very difficult to evaluate
  with constraints with the same characteristics. In these cases, common sense indicates that one
  should try to minimize suitable progressive approximations with the hope that 
   successive partial minimizers would converge to the solution of the original problem. In many cases
  error bounds are not available, so that we know how to get closer to the true problem 
   but we cannot estimate distances between partial and final solutions. 

   The natural questions that arise are: With which precision we need to solve each partial 
  problem? How to choose the approximate  problem that should be addressed after 
   finishing each stage of the process?  For solving these questions one needs to consider 
   two different objectives: decreasing the objective function and increasing the precision. 
    It is natural to combine these objectives in a single merit function. 

  The papers  \cites{bkm2018,bkm2019,bkm2021,nkjmm} suggested that a good framework to address
   this problem is given by the Inexact Restoration approach of classical constrained optimization.
   The idea is that ``maximal evaluation precision" can be considered as a constraint of the problem
  depending on a precision variable $y$ that lies in an abstract set $Y$. The tools of 
  Inexact Restoration indicate an algorithmic path for modifying $y$ and decreasing the objective
   function in such a way that, hopefully, most iterations are performed with moderate precision
   and the overall computational cost is affordable.

   The present paper is the first contribution in which the Inexact Restoration framework is 
  applied to the case in which, not only the objective function but also the constraints are
  subject to uncertainty. An interesting feature of our approach is that  
   our method applied to the particular case 
   in which  exact evaluations are possible ($Y$ is a singleton, $g_f(y)=0$ and $g_h(y)=0$) 
    coincides with (a version of)  the classical Inexact Restoration method for smooth constrained
  optimization. Paradoxically, this nice feature motivates a challenging open problem: 
     Is it really necessary to use the IR approach both for the algebraic and the precision
  constraints?  From the aesthetic point of view our ``double IR" strategy seems to be 
   attractive but it cannot be discarded that using different underlying strategies for the
  algebraic constraints could result in more efficient algorithms. 
     
 An important possible branch of application of the theory of this work is Stochastic Optimization. In \cites{bcrz,cor,crz}  only the objective function is stochastic whereas the constraints are deterministic.   However, some nontrivial adaptations of the main algorithm may be necessary to consider the   specific contribution of stochasticity.  The paper \cite{bkmr} presents several successful   applications of the IR approach and, in particular, show the way in which functions $h(y)$ and merit functions can be  defined for that type of problems. The application  to noisy derivative free optimization \cite{sxxn}, on the other hand, will be also the subject of    future research.
  \\           

\noindent
{\bf Acknowledgments}
 
The authors are deeply indebted to two anonymous reviewers, whose comments and recommendations helped  a lot
  to improve the readability of this paper.

  \appendix

\newcounter{appendix} 
\stepcounter{section}
 \setcounter{appendix}{0}
  \stepcounter{appendix}
\renewcommand{\thesection}{\Alph{section}}
  
\section*[A]{Appendix A: Proofs of auxiliary results} \label{apend}

\subsection{Proof of Lemma~\ref{restorationcontrolc}}

\begin{proof}
	Using \eqref{ctaylor} for
	$x_2=z^{trial}$ and $x_1=z^\ell$ we have that
	\[
	c(z^{trial},w^{i+1}) \leq c(z^\ell,w^{i+1}) + \nabla_x c(z^\ell,w^{i+1})^T (z^{trial} - z ^\ell)
	+ L_c \| z^{trial} -z ^\ell \| ^ 2.
	\]
	Then, taking
	\[
	v = \frac{1}{2} (z^{trial} - z^\ell)^T B_\ell (z^{trial} - z^\ell) + \frac{\sigma}{2} 
	\|z^{trial} - z^\ell\|^2 ,
	\]
	we obtain:
	\begin{equation*}\begin{split}
			c(z^{trial},w^{i+1}) - c(z^\ell,w^{i+1})   &\leq  \nabla_x c(z^\ell,w^{i+1})^T (z^{trial} - z^\ell) + L_c \|z^{trial} - z^\ell\|^2 \\
			&= \nabla_x c(z^\ell,w^{i+1})^T (z^{trial} - z^\ell) + v-v\\&\quad + L_c \|z^{trial} - z^\ell\|^2 \\
			& =  \nabla_x c(z^\ell,w^{i+1})^T (z^{trial} -  z^\ell)  + v\\&\quad   + 
			\left( L_c  -  \frac{\sigma}{2}  \right) \|z^{trial} - z^\ell\|^2\\
			&\quad - \frac{1}{2} (z^{trial} - z^\ell)^T B_\ell (z^{trial} - z^\ell).\end{split}
	\end{equation*}

	Since  $\|B_\ell\| \leq M $, we have that
	$|(z^{trial} - z^\ell)^T B_\ell (z^{trial} - z^\ell)| \le M \|z^{trial} - z^\ell\|^2$, so:
	\begin{equation}\label{deltac}
 \begin{split}
     c(z^{trial},w^{i+1}) - c(z^\ell,w^{i+1})   &\le  \left[ \nabla_x c(z^\ell,w^{i+1})^T (z^{trial} -  z^\ell)  + v \right] \\
  & \quad + 
		\left( L_c +\frac{M}{2}-  \frac{\sigma}{2} \right) \|z^{trial} - z^\ell\|^2.
 \end{split}
	\end{equation}

	By  \eqref{eqreduzc}, 
	\begin{equation}\label{deltac2}
		c(z^{trial},w^{i+1}) - c(z^\ell,w^{i+1})  \le \left( L_c +\frac{M}{2}-  \frac{\sigma}{2} \right) \|z^{trial} - z^\ell\|^2.
	\end{equation}
	
	Therefore, 
	we obtain that, if $ {\sigma} \geq  \bar{\sigma}$, \eqref{desc} is fulfilled.   \end{proof}   

\subsection{Proof of Lemma~\ref{thefo}}

\begin{proof} Define 
	$$ u=z^{\ell+1}-[\nabla_x c(z^\ell,w^{i+1})+B_\ell (z^{\ell+1}-z^\ell)] \,\, \text{ and  } \,\, w=u-\sigma (z^{\ell+1} - z^\ell). $$
	By   \eqref{projnula}, 
	\begin{equation}\label{now}
 \begin{split}
     \norm{ P_\Omega (u) - z^{\ell+1} } &=\norm{ P_\Omega (u)- P_\Omega (w)+ P_\Omega (w) - z^{\ell+1} } \\& \leq \norm{ P_\Omega (u)   - P_\Omega ( w )}+ \kappa_R \|z^{\ell+1} - z^\ell\|.
 \end{split}
	\end{equation}
	By the non-expansivity projections, we have that
	\begin{equation*} 
		\norm{  P_\Omega (u)   - P_\Omega ( w )} \leq \norm{ u - w } = \sigma \norm{ z^{\ell+1} - z^\ell }.
	\end{equation*} 
	So, \eqref{now} implies that
	\begin{equation}\label{now1}
		\norm{  P_\Omega (u) - z^{\ell+1} }  \leq  (\sigma+ \kappa_R) \norm{ z^{\ell+1} - z^\ell }.
	\end{equation} 
	Now, define  $v= z^{\ell+1} - \nabla_x c (z^{\ell+1},w^{i+1})$. Using  \eqref{now1} and, again, the non-expansivity of projections,
	we obtain: 
	\begin{equation}\label{now2}
		\begin{split}\norm{ P_\Omega (v ) - z^{\ell+1} }  & \leq \norm{ P_\Omega (v) - P_\Omega ( u)  } + \norm{ P_\Omega ( u) - z^{\ell+1} }  \\ 
			& \leq \norm{ v-u } + (\sigma+\kappa_R) \norm{ z^{\ell+1} - z^\ell }\\
			& \leq \norm{ \nabla_x c(z^\ell,w^{i+1})- \nabla_x  c(z^{\ell+1},w^{i+1}) + B_\ell (z^{\ell+1} - z^\ell) } \\
   & \quad + (\sigma+\kappa_R) \norm{ z^{\ell+1} - z^\ell }.
		\end{split}
	\end{equation}

	By   \eqref{lipsgradc}, we have that
	$\norm{ \nabla_x c(z^{\ell+1},w^{i+1}) - \nabla_x c(z^\ell,w^{i+1})} \leq L_c \norm{z^{\ell+1} - z^\ell}.$ So, since $  \norm{ B_\ell } \leq M $, 
	\begin{equation} \label{now3}
		\begin{split}
			\norm{ P_\Omega (v) - z^{\ell+1}  }  &  \leq \norm{ \nabla_x c(z^\ell,w^{i+1})- \nabla_x  c(z^{\ell+1},w^{i+1})}  +\norm{  B_\ell (z^{\ell+1} - z^\ell)}\\
   &\quad+ (\sigma+\kappa_R) \norm{ z^{\ell+1} - z^\ell }\\
			& \leq L_c \norm{z^{\ell+1} - z^\ell} + M \norm{z^{\ell+1} - z^\ell}+ (\sigma+\kappa_R) \norm{ z^{\ell+1} - z^\ell } \\ & = (L_c+M +\sigma+\kappa_R) \norm{z^{\ell+1} - z^\ell}. \end{split}
	\end{equation} 
	Therefore, 
	recalling that, by Lemma 
	\ref{restorationcontrolc},
	we have that  $\sigma \leq \max\{10\bar{\sigma}, \sigma_{max}\}$, we deduce  \eqref{graleq0}, as desired.
\end{proof}

\subsection{Proof of Lemma \ref{resumoGencan} }

\begin{proof} If $\|h(x^k, w^{i+1})\|=0$, then $z^0=x^k$ satisfies \eqref{restctargetz} and \eqref{restnablacz}. Assume now that $\|h(x^k, w^{i+1})\|>0$ and \eqref{restnablacz} is not true for the first $\ell$ iterations
	of Step~4 of RESTA. Then, for all $j \in \{0, 1, \ldots, \ell\}$ we have that
	\begin{equation} \label{dentro1} 
		\norm{ P_{\Omega}(z^j  -  \nabla_x c(z^j,w^{i+1}) )-z^j } > \epsilon_c. 
	\end{equation}

	By  \eqref{graleq0}, for all $j \in \{0, 1, \ldots, \ell-1\}$, 
	\begin{equation} \label{dentro2} 
		\norm{ P_\Omega \left( z^{j+1} - \nabla_x c (z^{j+1},w^{i+1}) \right) - z^{j+1} } \leq  c_{P_{\Omega}}   \norm{ z^{j+1}-z^j}.  
	\end{equation} 
	Then, by (\ref{dentro1}) and (\ref{dentro2}),  
	\begin{equation}\label{descintermedio}
\begin{split}
    \ell\epsilon_c^2 &=  \;\sum_{j=0}^{\ell-1} \epsilon_c^2 \le   \;\sum_{j=0}^{\ell-1} \norm{ P_{\Omega}(z^{j+1}  -  \nabla_x c(z^{j+1},w^{i+1})) - z^{j+1} }^2\\&  \leq  c_{P_{\Omega}}^2    \sum_{j=0}^{\ell-1} \norm{z^{j+1}-z^{j} }^2.
\end{split}		
	\end{equation} 
	
	On the other hand,  
	\begin{equation*} 
		c(z^{\ell},w^{i+1})    =  { \sum_{j=0 }^{\ell-1} }\left[  c(z^{j+1},w^{i+1}) -c(z^j,w^{i+1})  \right] + c(z^0, w^{i+1}). 
	\end{equation*}
	
	By  \eqref{desc} at Step~5.2 of Algorithm RESTA, we have that
	$c(z^{j+1},w^{i+1}) \leq c(z^j,w^{i+1}) - \alpha_R \|z^{j+1}- z^j\|^2$, 
	for all  $j \in \{0, 1, \ldots, \ell-1\}$. Therefore,  
	\begin{equation*}
 \begin{split}
		c(z^{\ell},w^{i+1}) - c(z^0, w^{i+1})  &= 
		{\sum_{j=0 }^{\ell-1} }\left[  c(z^{j+1},w^{i+1}) -c(z^j,w^{i+1}) \right] \\&\le - \alpha_R  {\sum_{j=0 }^{\ell-1} } \|z^{j+1}- z^j\|^2 .
  \end{split}
	\end{equation*}
	By \eqref{descintermedio} and the fact that, at Step~3 of 
	RESTA, we choose $z^0$ such that $c(z^0, w^{i+1}) \leq c(x^k, w^{i+1})$, we have that
	\begin{equation*}
		c(z^{\ell}, w^{i+1}) \leq c(z^0, w^{i+1})   -  \alpha_R \frac{\ell \epsilon_c^2} { c_{P_{\Omega}}^2  }   \leq c(x^k, w^{i+1})   -  \alpha_R   \frac{\ell \epsilon_c^2} { c_{P_{\Omega}}^2  }.
	\end{equation*}
	Therefore, if 
	\begin{equation}\label{hiplrestoration}
		c(x^k, w^{i+1})   - \frac{\alpha_R\epsilon_c^2  }{c_{P_\Omega}^2}  \ell \leq c_{target}
	\end{equation} 
	we would have that $c(z^{\ell}, w^{i+1})\leq c_{target}$ and the stopping condition \eqref{thesto} at Step 4would be fulfilled.
	
	Moreover, \eqref{hiplrestoration} occurs if and only if  
	\begin{equation}\label{boundlrestoration}
		\frac{c_{P_\Omega}^2}{\alpha_R\epsilon_c^2  } \left[c(x^k, w^{i+1})   -c_{target} \right]   \leq  \ell.\end{equation}
	Using the definitions of
	$\epsilon_c$ and $c_{target}$, given in \eqref{parafeas}, we obtain that
	\eqref{boundlrestoration} is equivalent to
	\begin{equation*}
		\ell \ge  \frac{c_{P_\Omega}^2}{\alpha_R \left[r_{feas}\|h(x^k, w^{i+1})\| \right]^2  } \left[c(x^k,w^{i+1})   -r^2 c(x^k,w^{i+1}) \right] = \frac{c_{P_\Omega}^2(1-r^2)}{2\alpha_R\,  r_{feas}^2} .   
	\end{equation*}
	Therefore, if 
	$\ell \ge \frac{c_{P_\Omega}^2(1-r^2)}{2\alpha_R\,  r_{feas}^2}  $ and
	\eqref{restnablacz} has not been fulfilled before, we have that
	$c(z^{\ell},w^{i+1})   \leq c_{target}$. Therefore, 
	we have that in at most  $C_{rest}$ sub-iterations of RESTA
	either \eqref{restctargetz} holds or  \eqref{restnablacz}  would have been obtained before.
	This completes the proof.     \end{proof}    

\subsection{Proof of Lemma \ref{lemascurto}}

\begin{proof} If $z^{\ell+1} = z^\ell$, (\ref{cotazeta0}) is trivial. If
	$z^{\ell+1} \neq z^\ell$, define
	$		v = \frac{ z^{\ell+1} - z^\ell  }{\norm{   z^{\ell+1} - z^\ell  } }.$
	
	Since $z^{\ell+1}$ is an approximate minimizer of  \eqref{subproblemfeas},
	then, by Assumption~\ref{reduzc}, 
	\begin{equation*}
		\nabla_x c(z^\ell,w^{i+1})^T (z^{\ell+1} - z^\ell) \leq -\frac{1}{2} (z^{\ell+1}- z^\ell)^T (B_\ell + \sigma I) (z^{\ell+1}- z^\ell)  .
	\end{equation*}  
	By Step~5 of RESTA, we have that
	$\sigma \geq \sigma_{min}$ and $B_\ell + \sigma_{min} I$ is positive definite, so 
	$\nabla_x c(z^\ell,w^{i+1})^T v<0$.
	
	Now, consider the function 
	$\varphi: \mathbb{R}_+ \to  \mathbb{R} $ defined by
	\begin{equation}\label{defvarphi}
		\varphi(t) = t  \nabla_x c(z^\ell,w^{i+1})^T v + \frac{t^2}{2}  v^T (B_\ell + \sigma I)v. 
	\end{equation}
	The unconstrained minimizer of $\varphi(t)$ is
	\begin{equation*}
		t^* = - \frac{ \nabla_x  c(z^\ell,w^{i+1})^T v   }   {v^T (B_\ell + \sigma I)  v}    \leq   - \frac{ \nabla_x  c(z^\ell,w^{i+1})^T v   }   {v^T (B_\ell + \sigma_{min} I)  v} \leq \frac{  \norm{ \nabla_x   c(z^\ell,w^{i+1}) } \|v\|  }{\lambda_1(B_\ell+ \sigma_{min} I) \|v\|^2 },
	\end{equation*}
	where $\lambda_1(B_\ell+ \sigma_{min} I) > 0$ is the smaller 
	eigenvalue of  $B_\ell+ \sigma_{min} I$.
	As  $ \| v \| = 1$,  by Step~5 of RESTA,  
	\begin{equation} \label{testrella}
		t^* \leq \norm{ (B_\ell + \sigma_{min} I)^{-1}} \norm{ \nabla_x c(z^\ell,w^{i+1})} \leq  M   \norm{ \nabla_x c(z^\ell,w^{i+1})}. 
	\end{equation}
	
	Let  $\bar{t}$ be the minimizer of  $\varphi(t)$ subject to
	$z^{\ell} + t v \in \Omega$. By the convexity of $\Omega$ we have that
	$\bar{t} \leq t^*$. Moreover, by construction, 
	$ z^\ell + \bar{t} v = z^{\ell+1}_*$. So, by  $\|v\|=1$, Assumption~\ref{assumpdistminunidimen},   
	and \eqref{testrella}, we have that 
	\begin{equation*}
		\norm{ z^{\ell+1} - z^\ell} \leq \kappa_{\varphi}	\norm{ z^{\ell+1}_* - z^\ell}  = \kappa_{\varphi} \bar{t} \leq \kappa_{\varphi} t^*\leq \kappa_{\varphi} M   \norm{ \nabla_x   c(z^\ell,w^{i+1})}.
	\end{equation*} 
	Now, by \eqref{boundgrau}, 
	$\|  \nabla_x h (z^\ell,w^{i+1})^T \|=\|  \nabla_x h (z^\ell,w^{i+1}) \|  \leq C_{h}$, thus: 
	\begin{equation*} \label{cotazeta}
		\begin{array}{cl}  
			\norm{ z^{\ell+1} - z^\ell }& \leq  \kappa_{\varphi} M  \norm{\nabla_x h (z^\ell,w^{i+1})^T h(z^\ell,w^{i+1})} \\
			& \leq \kappa_{\varphi} M  \norm{\nabla_x h (z^\ell,w^{i+1})^T} \norm{h(z^\ell,w^{i+1})} \\
			& \le \kappa_{\varphi} M C_{h} \|h(z^\ell,w^{i+1})\|.  
		\end{array}
	\end{equation*} 
	
	By \eqref{desc} we have that 
	$\|h(z^{\ell+1},w^{i+1})\| \leq \|h(z^\ell,w^{i+1})\|$  and, 
	by the choice of  $z^0$ at Step~3 of RESTA, we have that
	$\|h(z^{0},w^{i+1})\| \leq \|h(x^k,w^{i+1})\|$. Then,
	$\|h(z^\ell,w^{i+1})\| \le \|h(x^k,w^{i+1})\|$ for all $\ell$. Therefore, 
	\eqref{cotazeta0} holds.
	
\end{proof}

\subsection{Proof of Lemma \ref{Navhrest}}

\begin{proof}  
	For a fixed $i$, by Lemma  \ref{resumoGencan}, after at most
	$C_{rest}$ steps we find $z^{\ell}$ satisfying \eqref{restctargetz} or
	\eqref{restnablacz}. At each of these steps,  $\nabla_x c $ is evaluated
	only once, therefore, the total number of evaluations
	of  $\nabla_x c $ is bounded by $C_{rest}.$
	
	By Lemma \ref{restorationcontrolc}, 
	$\sigma \ge \bar{\sigma} $ implies that
	\eqref{desc}  is fulfilled and, so,  $z^{\ell+1}$ is well defined.
	By  \eqref{updatemufeas},  $\sigma$ is increased according to
	$\sigma \in [2\sigma, 10\sigma]$. Therefore, as the initial value
	of $\sigma$ is not smaller than  $\sigma_{min}$, we have that
	after   ${n}_{\sigma}$ trials we will have that 
	$  \sigma \geq 2^{n_\sigma} \sigma_{\min}$. Therefore, 
	we have that 
	$\sigma \geq 2^{n_\sigma} \sigma_{\min} \geq \bar{\sigma}$ and, so, $z^{\ell+1}$ is obtained.   
	
	Therefore, the descent condition 
	\eqref{desc} is tested at most  ${n}_\sigma$ times for each value of
	$\ell$. Consequently, $h$ is evaluated at most ${n}_\sigma$ times for every  $\ell$.
	So, the condition 	
	\eqref{desc} is tested at most  $C_{rest}n_\sigma$ times for all fixed $w^{i+1}$.
	
	Finally, observe that, by Step~2, after at most 
	$N_{prec}$ trials we have that $g_h(w^{i+1}) \leq\bar{\epsilon}_{prec}$. 
	In this case, the process would ensure that conditions \eqref{thesto} or \eqref{thesto2} hold   and, so, the Restoration Phase 
	would be finished. Moreover, only one additional evaluation of $h$ is performed at each
	update of $w^{i+1}$. Then,
	we obtain the desired result.    
\end{proof}            

\subsection{Proof of Lemma \ref{valecondbetasobreeps}} 

\begin{proof} If $(x_R^k,y_R^k)$ is computed by a problem-dependent procedure, the result is true by \eqref{beta3}. Now, let us consider that $(x_R^k,y_R^k)$ is computed by RESTA. For given $k$, let 
	$i$ be such that  $w^{i+1}=y^k_R$ and let $N_{Rk}$ be 
	the number of sub-iterations performed for the minimization of
	$c(z,w^{i+1})$. By Lemma  \ref{Navhrest} we have that  
	$N_{Rk} \leq N_{RESTA}$.
	Then, by Lemma~\ref{lemascurto} and the choice of $z^0$ at Step~3 of RESTA, we have that
	\begin{equation*} \begin{array}{cl}
			\norm{ x_R^k-x^k}  & \leq \|z^0-x^k\|+{\sum_{l=1}^{N_{Rk}}} \norm{  z^\ell-z^{\ell-1} }\\
			& \leq \beta_c  \|h(x^k,y^k_R)\| +{\sum_{l=1}^{N_{Rk}}}  C_s   \|h(x^k,y^k_R)\|  \\ & 
			\leq \beta_c  \|h(x^k,y^k_R)\| +N_{RESTA}  C_s   \|h(x^k,y^k_R)\|. \end{array}
	\end{equation*} 
	Therefore, 
	we obtain  the desired result. 
\end{proof}

\subsection{Proof of Lemma  \ref{boudfbyinfeas}}

\begin{proof} By \eqref{lipsf} we have that 
	$  |f(x_R^k,y^k_R) - f(x^k,y^k_R)|  \leq  L_f \norm{ x_R^k - x^k }$.
	Then, by \eqref{restorationf},
	\begin{equation*} \begin{array}{cl}
			f(x_R^k,y^k_R) - f(x^k,y^k)  &  \le   f(x_R^k,y^k_R) -  f(x^k,y^k_R)   +     f(x^k,y^k_R) - f(x^k,y^k)  \\ & \le  L_f \norm{ x_R^k - x^k } + \beta g(y^k).\end{array}
	\end{equation*}
	By Lemma \ref{valecondbetasobreeps} we have that 
	$\norm{ x_R^k-x^k}  \le \beta_R  \|h(x^k,y^k_R)\| $. Then,  
	\begin{equation*}\begin{array}{cl}
			f(x_R^k,y^k_R) - f(x^k,y^k)   & \le  L_f \beta_R \|h(x^k,y^k_R)\|  + \beta g(y^k)\\ & \leq [ L_f\beta_R+ \beta] \left[ \|h(x^k,y^k_R)\|+g(y^k)  \right] .\end{array}
	\end{equation*}
	Thus, we have the desired result. 
\end{proof}

\subsection{Proof of Lemma \ref{tetabemdef}}

\begin{proof} 
	At each iteration $k$ of BIRA we have two options, according
	to the fulfillment of  \eqref{checkpenalty}. If (\ref{checkpenalty}) holds, 
	we define 
	$\theta_{k+1} = \theta_k $, therefore $\theta_{k+1}$ is well defined and 
	does not increase with respect to $\theta_k$. Moreover, in this case \eqref{penaltetamas} is equivalent to \eqref{checkpenalty}, so it is fulfilled.
	
	In the second case,   $\theta_{k+1}$ is defined by
	\eqref{updatepenaltyparameter} at Step~2, according to:
	{ \begin{equation*} 
 \begin{array}{c}
      			{\theta}_{k+1} = \frac{(1+r) [  \|h(x^k,y^k_R)\| -  \|h(x^k_R,y_R^k)\|   + g(y^k) - g(y_R^k)]} { 2 \left[ f(x_R^k,y_R^k) - f(x^k,y^k_R) + \|h(x^k,y^k_R)\| - \|h(x^k_R,y_R^k)\| + g(y^k) - g(y_R^k)\right]}.
 \end{array}
	\end{equation*} }
	
	Let us show that both the numerator and the denominator of this expression are positive
	and that the quotient
	is smaller than $\theta_k$. 
	
	By the restoration step and the assumptions G1--G12, we have that 
	$g(y_R^k) \leq r g(y^k)$, so $g(y_R^k) - g(y^k) \leq  0$.
	Therefore, as  $\frac{1-r}{2} \in (0,1)$, we have that
	\begin{equation} \label{rg}
		g(y_R^k) - g(y^k) \leq \frac{1-r}{2} [g(y_R^k) - g(y^k) ].
	\end{equation}
	Moreover, if the execution of BIRA is not stopped declaring Restoration Failure, the restoration always guarantees that 
	$ \|h(x^k_R,y_R^k)\| \le r \|h(x^k,y^k_R)\|$. Therefore,   
	\begin{equation} \label{rh}
		\|h(x^k_R,y_R^k)\| - \|h(x^k,y^k_R)\| \le \frac{1-r}{2} \left[  \|h(x^k_R,y_R^k)\| - \|h(x^k,y^k_R)\|\right]. 
	\end{equation} 
	Now, the equalities in 
	\eqref{rg}  and \eqref{rh} only take place if
	$\|h(x^k_R,y_R^k)\| = \|h(x^k,y^k_R)\|=g(y_R^k) = g(y^k)=0$.
	In this case, if $(x^k_R,y_R^k)$ is computed by the PDP, by \eqref{beta3}, we have that $(x^k_R,y_R^k)=(x^k,y^k)$. On the other hand, if RESTA is used, since $g(y^k)=0$, by Step~2, we would have that
	$w^i=y^k$  for all $i$, implying that  $y^k_R=y^k$. So,
	$\|h(x^k,y^k)\|=\|h(x^k,y^k_R)\|=0$ and, by Step~1 of RESTA, we also have that
	$(x^k_R,y_R^k)=(x^k,y^k)$. 
	In this case, \eqref{checkpenalty} would be trivially fulfilled  and we would have that
	$\theta_{k+1} = \theta_k$. Then, at least one of the conditions
	\eqref{rg} or \eqref{rh} is strictly satisfied. 
	This proves that, when $\theta_{k+1} \neq \theta_k$, the numerator of  \eqref{updatepenaltyparameter} is positive.

	Now let us analyze the expression $\Phi(x_R^k, y_R^k,  \theta)-\Phi(x^k, y^k_R,  \theta)$ as a function of $\theta$. By the definition of he merit function in \eqref{defmeritfunction}, we have that 
	\begin{equation} \label{equivpenal}
		\begin{split}
			\Phi(x_R^k, y_R^k,  \theta)-\Phi(x^k, y^k_R,  \theta)&=
			\theta \left[ f(x_R^k,y_R^k) - f(x^k,y^k_R) +\|h(x^k,y^k_R)\|\right.\\ &\quad  \left. - \|h(x^k_R,y_R^k)\| \right]- [  \|h(x^k,y^k_R)\| -  \|h(x^k_R,y_R^k)\|   ],
		\end{split} 
	\end{equation}
which is linear with respect to $\theta$ and its slope  is less or equal half the denominator of \eqref{updatepenaltyparameter}.  Moreover, this slope must be positive, otherwise, by \eqref{errata}, \eqref{checkpenalty} would hold for all non-negative $\theta$. So the expression of $\theta_{k+1}$ is well defined and  $\Phi(x_R^k, y_R^k,  \theta)-\Phi(x^k, y^k_R,  \theta)$ is a increasing bijection from $\R$ to $\R$.
	
	When $\theta =0 $ we have that 
	\begin{equation*}
		\Phi(x_R^k, y_R^k,  0)-\Phi(x^k, y^k_R,0)  =   \left[ \|h(x^k_R,y_R^k)\| - \|h(x^k,y^k_R)\| \right].
	\end{equation*} 
	Since one of the inequalities in 
	\eqref{rg} or  \eqref{rh} is strict, by \eqref{errata}, we have that
	\begin{equation} \label{penalty0}
 \begin{array}{c}
\Phi(x_R^k, y_R^k,  0)-\Phi(x^k, y^k_R,0)\\   < \frac{1-r}{2} \left[  \|h(x^k_R,y_R^k)\| - \|h(x^k,y^k_R)\|  + g(y_R^k) - g(y^k) \right].
 \end{array}
	\end{equation} 
	However, if 
	\eqref{checkpenalty} does not hold, we have that
	\begin{equation} \label{notpenalty}
 \begin{array}{c}
     \Phi(x_R^k, y_R^k,   \theta_k)-\Phi(x^k, y^k_R,  \theta_k)\\ >   \frac{1-r}{2}\left[ \|h(x^k_R,y_R^k)\|  - \|h(x^k,y^k_R)\| + g(y_R^k)-g(y^k) \right].
 \end{array}
	\end{equation}  
	So there exists only one value of 
	$\theta \in (0, \theta_k)$ verifying   
 \begin{equation*}
 \begin{array}{c}
     \Phi(x_R^k, y_R^k,  \theta)-\Phi(x^k, y^k_R,  \theta)\\  = \frac{1-r}{2} \left[  \|h(x^k_R,y_R^k)\| - \|h(x^k,y^k_R)\|  + g(y_R^k) - g(y^k) \right].
 \end{array}
 \end{equation*}
By \eqref{equivpenal}, this value of $\theta$ is greater than or equal to 
	$\theta_{k+1}$ computed in \eqref{updatepenaltyparameter}. Therefore we also have that \eqref{penaltetamas} holds.
	So, the
	proof is complete.         
\end{proof}

\subsection{Proof of Lemma \ref{lemabarteta}}

\begin{proof} 
	It is enough to prove that 
	$\theta_{k+1}$ is bounded below by $\bar{\theta}$ when
	it is defined by \eqref{updatepenaltyparameter}.
	
	Equivalently, we need to show that 
	$\frac{1}{\theta_{k+1}}$ is bounded above in this situation. In fact,
	\begin{equation}\label{thetaawayzero0}
		\begin{array}{cl}
			\frac{1}{{\theta}_{k+1}} & =
			\frac  { 2 [ f(x_R^k,y_R^k) - f(x^k,y^k_R)+  \|h(x^k,y^k_R)\| - \|h(x^k_R,y_R^k)\| + g(y^k)- g(y_R^k)]}
			{(1+r) [ \|h(x^k,y^k_R)\| - \|h(x^k_R,y_R^k)\| +g(y^k) -g(y_R^k)]} \\
			&= \frac{2}{1+r} \left[ \frac{f(x_R^k,y_R^k) - f(x^k,y^k_R)}{ \|h(x^k,y^k_R)\| - \|h(x^k_R,y_R^k)\| + g(y^k) -g(y_R^k)} + 1\right].
		\end{array}
	\end{equation}

	By Step~1 of RESTA or \eqref{erre3} when using a PDP, 
	$$ - \|h(x^k_R,y_R^k)\|-g(y_R^k) \geq - r\|h(x^k,y_R^k)\|-rg(y^k),$$ therefore
	$$\begin{array}{c}
		\|h(x^k,y_R^k)\| + g(y^k) - \|h(x^k_R,y_R^k)\|-g(y_R^k) \\
  \geq \|h(x^k,y_R^k)\| + g(y^k) - r\|h(x^k,y_R^k)\|-rg(y^k)\\ = (1-r)(\|h(x^k,y_R^k)\| + g(y^k))>0.
	\end{array}$$
	Positivity necessarily takes place, otherwise we would have that 
	$(x^k_R,y_R^k)=(x^k,y^k)$ and 
	$\theta_{k+1}=\theta_k$. Thus,   
	\begin{equation*}
	    \begin{split} 0&<\frac{1}{\|h(x^k,y_R^k)\| + g(y^k) - \|h(x^k_R,y_R^k)\|-g(y_R^k)} \\ &\leq \frac{1}{ (1-r)(\|h(x^k,y_R^k)\| + g(y^k))}. \end{split}
	\end{equation*}

	On the other hand, by 
	\eqref{lipsf}, we have that
	$ |f(x_R^k,y_R^k) - f(x^k,y_R^k)|  \leq  L_f \|x_R^k - x^k\|$. Then, by 
	\eqref{thetaawayzero0}, 
	\begin{equation*}
		\begin{array}{cl}
			\frac{1}{{\theta}_{k+1}} 
			& \le \frac{2}{1+r} \left[ \frac{L_f \|x_R^k-x^k\| }{ (1-r)(\|h(x^k_R,y_R^k)\| + g(y^k))} + 1 \right].  
		\end{array}
	\end{equation*}  
	
	By \eqref{cotaxyfeas} in Theorem \ref{teofinres}, 
	there exists a positive constant  $\beta_R= O(1)$ such that
	$\|x_R^k - x^k\|  \le \beta_R \|h(x^k,y^k_R)\|  $.    Then, since $g(y^k) \geq 0$,
	\begin{equation*}
 \begin{split}
			\frac{1}{\theta_{k+1}} & \leq   \frac{2}{1+r} \left[ \frac{L_f \beta_R \|h(x^k,y^k_R)\| }{(1-r)(\|h(x^k,y_R^k)\| + g(y^k))} + 1 \right] \\
   & \leq  \frac{2}{1+r} \left[ \frac{L_f \beta_R (\|h(x^k,y^k_R)\| + g(y^k)) }{(1-r)(\|h(x^k,y_R^k)\| + g(y^k))} + 1 \right]= \frac{2}{1+r} \left[ \frac{L_f \beta_R }{1-r} + 1 \right]. 
   \end{split} 
	\end{equation*}
	
	The inequality above implies that, when $\theta_k$ is updated,  
	$\{ \frac{1}{\theta_k}\}$ is bounded, so  $\{\theta_k\}$ is bounded away from zero, with 
	$\bar{\theta} > 0$ as lower bound.   
\end{proof}

\subsection{Proof of Lemma \ref{lemmaNreg0}}

\begin{proof} 
	By \eqref{ftaylor}, we have that 
	$f(x,y^{k+1}) \leq  f(x_R^k,y^{k+1}) +\nabla_x f(x_R^k,y^{k+1})^T (x - x_R^k) + L_f \|x - x_R^k\|^2$.

	Then, since
	$\|H_k\| \leq M$, we have that  
	{ \begin{equation*}
			\begin{split}
				f(x,y^{k+1}) & \leq    f(x_R^k,y^{k+1}) + \nabla_x f(x_R^k,y^{k+1})^T (x - x_R^k) \\
    & \quad +  \frac{1}{2} (x - x_R^k)^T H_k (x - x_R^k)  - \frac{1}{2} (x - x_R^k)^T H_k (x - x_R^k)  \\				%
				&  \quad  + \tilde{\alpha}\|x - x_R^k\|^2 - \tilde{\alpha}\|x - x_R^k\|^2  + L_f \|x - x_R^k\|^2\\
				&  \leq  f(x_R^k,y^{k+1}) +  \nabla_x f(x_R^k,y^{k+1})^T (x - x_R^k) \\
    &\quad + \frac{1}{2} (x - x_R^k)^T H_k (x - x_R^k)+ M \|x - x_R^k\|^2\\
				&  \quad   + \tilde{\alpha}\|x - x_R^k\|^2 - \tilde{\alpha}\|x - x_R^k\|^2 + L_f \|x - x_R^k\|^2\\
				& \leq   f(x_R^k,y^{k+1}) +  \nabla_x f(x_R^k,y^{k+1})^T (x - x_R^k)  \\
    & \quad  + \frac{1}{2} (x - x_R^k)^T H_k (x - x_R^k)+ (M + \tilde{\alpha}+L_f) \|x - x_R^k\|^2 \\
				& \quad - \tilde{\alpha}\|x - x_R^k\|^2. 
    \end{split}
	\end{equation*} }

	Taking $\mu \geq C_\mu$, by Assumption~\ref{bajasub}, 
	\begin{equation*} \begin{array}{cl}
			f(x,y^{k+1}) &   \leq   f(x_R^k,y^{k+1}) - \tilde{\alpha}\|x - x_R^k\|^2 + \left[ \nabla_x f(x_R^k,y^{k+1})^T (x - x_R^k) \right.\\ 
    & \left. \quad + \frac{1}{2} (x - x_R^k)^T H_k (x - x_R^k) + \mu \|x - x_R^k\|^2 \right] \\
				& \leq   f(x_R^k,y^{k+1})  - \tilde{\alpha}\|x - x_R^k\|^2.\end{array}
	\end{equation*} 
	
	Since 
	$\alpha \leq \tilde{\alpha}$ and $y^{k+1}=y_R^k$,
	\eqref{desfinte} necessarily holds. Moreover, by  (\ref{defalfatil}), we have that 	
	\begin{equation}\label{desmu}
		f(x,y^{k+1})-f(x_R^k,y^{k+1}) \leq - \frac{ 1 - \bar{\theta}}{\bar{\theta}} (\kappa_T+L_h) \|x - x_R^k\|^2.
	\end{equation}

	Let us prove that  
	\eqref{desmeri} also holds when  $\mu \geq C_\mu $. Note that 
	\begin{equation}\label{nowr} 
 \begin{array}{c}
			\Phi(x, y^{k+1}, \theta_{k+1}) - \Phi(x^k, y^{k+1},  \theta_{k+1}) \\
   =[\Phi(x, y^{k+1}, \theta_{k+1}) - \Phi(x_R^k, y^{k+1}, \theta_{k+1})]  \\ 
\quad   + [\Phi(x_R^k, y^{k+1},  \theta_{k+1}) - \Phi(x^k, y^{k+1},\theta_{k+1})]. 
   \end{array}  
	\end{equation} 

	Define   $v =  \Phi(x, y^{k+1},  \theta_{k+1}) - \Phi(x_R^k, y^{k+1},  \theta_{k+1}).$ 
	By the definition of $\Phi$ and \eqref{desmu}, the first term in the right-hand side of the equality above we have that 
	\begin{equation*}
		\begin{split} 
			v &  =  \theta_{k+1} [ f(x,y^{k+1}) - f(x_R^k,y^{k+1}) ] \\
   & \quad +  (1 - \theta_{k+1})\left[\|h(x,y^{k+1})\| - \|h(x_R^k,y^{k+1})\| \right] \\
			& \leq \theta_{k+1}\left[ - \frac{ 1 - \bar{\theta}}{\bar{\theta}} (\kappa_T+L_h) \|x - x_R^k\|^2 \right]\\
   & \quad   +  (1 - \theta_{k+1}) \left[  \|h(x,y^{k+1})\|  -  \|h(x_R^k,y^{k+1})\|  \right].
		\end{split}
	\end{equation*} 
	By \eqref{utaylor} and \eqref{errotan},  
	\begin{equation*}
		\begin{split} 
			v &  \leq  \theta_{k+1} \left[ - \frac{ 1 - \bar{\theta}} {\bar{\theta}} (\kappa_T+L_h) \|x - x_R^k\|^2 \right]\\
&\quad +  (1 - \theta_{k+1})\left[\|\nabla_x h(x_R^k,y^{k+1})^T (x-x_R^k)\|  + L_h \|x-x_R^k\|^2\right]   \\
			&\leq \theta_{k+1} [ - \frac{ 1 - \bar{\theta}} {\bar{\theta}} (\kappa_T+L_h) \|x - x_R^k\|^2 ]\\&\quad    +  (1 - \theta_{k+1}) \left[ \kappa_T \|x-x_R^k\|^2+ L_h \norm{x - x_R^k}^2    \right].
		\end{split}
	\end{equation*} 
	Since  $\{\theta_k\}$  is bounded below by  $\bar{\theta}$, we have that
	$$ v  \leq \bar{\theta} \left[ - \frac{ 1 - \bar{\theta}} {\bar{\theta}} (\kappa_T+L_h) \|x - x_R^k\|^2 \right]  + (1 - \bar{\theta})  (\kappa_T+L_h)   \|x - x_R^k\|^2=0.$$ 
	
	Then, the first term in
	\eqref{nowr} is not positive. On the other hand, as 
	$y^{k+1}= y_R^k$,  \eqref{penaltetamas} is equivalent to
	$$\begin{array}{c}\Phi(x_R^k, y^{k+1},  \theta_{k+1}) - \Phi(x^k, y^{k+1},\theta_{k+1})\\ \leq  \frac{1-r}{2}\left[  \|h(x_R^k,y^{k+1})\| - \|h(x^k,y^{k+1})\|+ g(y^{k+1})- g(y^k) \right]. \end{array}$$
	Then, by \eqref{nowr}, 
	$$\begin{array}{c}\Phi(x, y^{k+1},  \theta_{k+1}) - \Phi(x^k, y^{k+1},\theta_{k+1}) \\ \leq  \frac{1-r}{2}\left[  \|h(x_R^k,y^{k+1})\| - \|h(x^k,y^{k+1})\|+ g(y^{k+1})- g(y^k) \right].\end{array}$$
	Therefore, if  $\mu \ge C_\mu$ , both \eqref{desfinte} and \eqref{desmeri}
	are fulfilled, guaranteeing that $x^{k+1}$ is well defined. 
\end{proof}

\subsection{Proof of Corollary \ref{lemmaNreg}}

\begin{proof}
	If  $x^{k+1}$ is computed  at the first $N_{acce}$ iterations of the Optimization Phase, we have that
	$\mu$ is increased at most $N_{accel}$ times, starting from a value limited above  by  $\mu_{max}$.  
	Therefore  $ \mu_k \le 10^{N_{acce}} \mu_{max}$.
	On the other hand, if $y^{k+1}=y^k_R$ and $\mu \ge C_\mu$, by Lemma \ref{lemmaNreg0}, the decrease conditions at Step~3.2 are satisfied. So, if the initial value of $\mu$ is greater than $C_\mu$, $\mu_k=\mu$.
	Otherwise, as $\mu_{new} \in [2 \mu, 10 \mu]$, we would have that
	$\mu_k \leq 10 C_\mu$.  Since $10^{N_{acce}} \mu_{max} \geq \mu_{max}$, we have that $\mu_k \leq \max\{10 C_{\mu},10^{N_{acce}} \mu_{max}\}$.

	Moreover, as the initial value of $\mu$ is not smaller than 
	$\mu_{min} $, after $N_{reg}$ updatings we have that $\mu \geq 2^{N_{reg}} \mu_{min}$. So,
	if  $y^{k+1}=y^k_R$ and  $2^{N_{reg}} \mu_{min} \geq C_\mu$ , or, equivalently,  $ N_{reg} + \log_2( \mu_{min}) \geq \log_2( C_\mu) $, 
	\eqref{desfinte} and \eqref{desmeri} are fulfilled. 
\end{proof}

\subsection{Proof of Lemma \ref{lemagradprojproporcionald}}

\begin{proof}
	Define $v = P_{D^{k+1}} (x^{k+1} - \nabla_x f(x_R^k,y^{k+1}) - H_k  (x^{k+1}-x_R^k) - 2 \mu_k  (x^{k+1}-x_R^k)     ) $.    
	
	By \eqref{projectionintermed},   
	\begin{equation} \label{intermedioprojintermed}
		\begin{array}{c}
			\| P_{D^{k+1}} (x^{k}_R - \nabla_x f(x^{k}_R,y^{k+1}) ) - x^{k}_R  \|  \\ =  \|   P_{D^{k+1}} (x^{k}_R - \nabla_x f(x^{k}_R,y^{k+1}) )  - v + v-  x^{k+1}+x^{k+1}- x^{k}_R  \| \\ 
			 \leq \|   P_{D^{k+1}} (x^{k}_R - \nabla_x f(x^{k}_R,y^{k+1}) )   - v \| + \| v -x^{k+1} \|+ \|x^{k+1}- x^{k}_R  \| \\
			 \leq \|   P_{D^{k+1}}  (x^{k}_R - \nabla_x f(x^{k}_R,y^{k+1}) )  - v \| + (\kappa+1) \| x^{k+1} -x_R^k \|.  \end{array}
	\end{equation} 
	By Step~3   of BIRA, we have that    
	$\| H_k\| \leq M$. 
	By the non-expansive property of projections, 
	\begin{equation}\label{eqprojauxiliar}
		\begin{array}{c}
			\|  P_{D^{k+1}}  (x^{k}_R - \nabla_x f(x^{k}_R,y^{k+1}) )  - v \|\\
    \leq \| x^{k}_R -x^{k+1}	+ H_k  (x^{k+1}-x_R^k)  + 2 \mu_k  (x^{k+1}-x_R^k)   \| \\ 
		 \leq \| H_k  (x^{k+1}-x_R^k)  \|  + (2 \mu_k+1)\|   (x^{k+1}-x_R^k)  \| \\ 
			 \leq  (  M + 2 \mu_k+1 ) \| x^{k+1}-x_R^k  \|. \end{array}
	\end{equation}
	
	By Corollary \ref{lemmaNreg},  
	we have that  $ \mu_k \le \bar{\mu} $. 
	Then, by  \eqref{intermedioprojintermed} and \eqref{eqprojauxiliar}: 
	\begin{equation*}
		\begin{array}{c}
			\| P_{D^{k+1}} (x^{k}_R - \nabla_x f(x^{k}_R,y^{k+1}) ) - x^{k}_R \|  \\ \leq \|   P_{D^{k+1}} (x^{k}_R - \nabla_x f(x^{k}_R,y^{k+1}) ) -v \|\ + (\kappa+1) \| x^{k+1} -x_R^k \|  \\
			\leq ( M + 2 \mu_k+1 ) \| x^{k+1}-x_R^k  \|  + (\kappa+1) \| x^{k+1} -x_R^k \| \\
			 =  ( M  +\kappa +2  \bar{\mu}+2) \| x^{k+1}-x_R^k  \|.
		\end{array}
	\end{equation*}
\end{proof}

\subsection{Proof of Lemma~\ref{lemasomagradprojproporcionald}}

\begin{proof}
	By Lemma \ref{lemagradprojproporcionald},  
	\begin{equation}\label{defcproj}
		\| P_{D^{j+1}} (x^{j}_R - \nabla_x f(x^{j}_R,y^{j+1}) ) -x^{j}_R\|   \leq  C_p \| x^{j+1}-x_R^j  \|, 
	\end{equation} 
	for all  $j$.
	
	Adding the first $k$ squared terms of \eqref{defcproj}, by  \eqref{somad}, we have that
	\begin{eqnarray*}
		\sum_{j=0}^{k}  \| P_{D^{j+1}}  (x^{j}_R - \nabla_x f(x^{j}_R,y^{j+1}) ) -x^{j}_R \|^2 & \leq &  \sum_{j=0}^k\ \left[  C_p \|   x^{j+1}-x_R^j      \| \right]^2 \\[-5pt]
		& = &  C_p^2 \sum_{j=0}^k\ \|   x^{j+1}-x_R^j  \|^2 \\
		& \leq & C_p^2 C_d.
	\end{eqnarray*}
	
\end{proof}            

\subsection{Proof of Lemma \ref{lemaNfeas}} 

\begin{proof}
	By \eqref{limsomah},   
	\begin{equation*}
		\sum_{j=0}^{k} [\|h(x^j,y^j_R)\|+ g(y^j)] \leq C_{feas}.
	\end{equation*}   
	Then, as  $0 \leq\|h(x^j_R,y^j_R)\| \leq r\|h(x^j,y^j_R)\|$ and $g(y^j) \geq 0$,
	\begin{equation*}
 \begin{split}
		rC_{feas} &\geq  \sum_{j=0}^{k} r\|h(x^j,y^j_R)\| + rg(y^j) \geq \; \;  { \sum_{\substack{j=0 \\\|h(x^j_R,y^j_R)\| >  \epsilon_{feas}}}^{k }  } \|h(x^j_R,y^j_R)\| \\& \geq N_{hinfeas} \epsilon_{feas}.
  \end{split}
	\end{equation*}
	So,  $  \floor{\frac{rC_{feas}}{\epsilon_{feas}}} \geq N_{hinfeas}$. Analogously, for 
	$g(y^j) > \epsilon_{prec}$, we have that
	\begin{equation*}
		C_{feas} \geq  \sum_{j=0}^{k}[ \|h(x^j,y^j_R)\| + g(y^j) ] \geq \; \;  { \sum_{\substack{j=0 \\ g(y^j) >  \epsilon_{prec}}}^{k }  }  g(y^j)  \geq N_{ginfeas} \epsilon_{prec},
	\end{equation*}
	therefore $  \floor{\frac{C_{feas}}{\epsilon_{prec}}} \geq N_{ginfeas}$.
	
	Finally, if
	$ \|h(x^k_R,y^k_R)\|>  \epsilon_{feas}$ or $ g(y^k_R)>  \epsilon_{prec}$ we have that 
	$\|h(x^j_R,y^j_R)\| + g(y^j_R) > \min\{\epsilon_{feas},\epsilon_{prec}\}$. Thus, as 
	$\|h(x^j_R,y^j_R)\| + g(y^j_R) \leq r(\|h(x^j,y^j_R)\|+ g(y^j))$,
		\begin{equation*}
  \begin{split}
			rC_{feas} &\geq  \sum_{j=0}^{k} r(\|h(x^j,y^j_R)\| + g(y^j))\\& \geq 
   { \sum_{\substack{j=0 \\\|h(x^j_R,y^j_R)\|+ g(y^j_R) > \min\{\epsilon_{feas},\epsilon_{prec}\}}}^{k }  } [\|h(x^j_R,y^j_R)\|+  g(y^j_R)]  \\&\geq N_{infeas} \min\{\epsilon_{feas},\epsilon_{prec}\},
   \end{split}
		\end{equation*}
	
	so $N_{infeas} \leq    \floor{\max\left\{\frac{rC_{feas} }{\epsilon_{feas}}, \frac{rC_{feas} }{\epsilon_{prec}}\right\}}.$
\end{proof}

\subsection{Proof of Lemma \ref{lemaNopt}}

\begin{proof}

	If during $N_{opt}$ iterations we have  $ \| P_{D^{j+1}} (x^{j}_R - \nabla_x f(x^{j}_R,y^{k+1}) ) -x^{j}_R \| > \epsilon_{opt}$, by \eqref{sumap2}, we have that  
	\begin{equation*}
		C_{proj} \geq \sum_{j = 0}^k  \| P_{D^{j+1}} (x^{j}_R - \nabla_x f(x^{j}_R,y^{k+1}) ) -x^{j}_R \|^2 \geq N_{opt} \epsilon_{opt}^2
	\end{equation*}
	Therefore, $  \floor{\dfrac{C_{proj}}{\epsilon_{opt}^2}} \geq N_{opt}$.
	
\end{proof} 

\bibliographystyle{amsplain}

\begin{thebibliography}{10}

\bibitem{ir7} R. Andreani, S. L. C. Castro, J. L. Chela,
  A. Friedlander, and S. A. Santos, \textit{An inexact-restoration method for
  nonlinear bilevel programming problems, }Comput. Optim. Appl. 43, pp. 307--328, 2009.
 

\bibitem{ahm} R. Andreani, G. Haeser, and J. M. Mart\'{\i}nez, 
   {\it On sequential
optimality conditions for smooth constrained optimization, }Optimization 60, pp. 627-641, 2011. 

 \bibitem{amrs} R. Andreani, J. M. Mart\'{\i}nez, A. Ramos, and P. J. S. Silva,  
   {\it Strict Constraint Qualifications and Sequential Optimality Conditions
for Constrained Optimization,} Math. Oper. Res. 43, pp. 693-717, 2018.
       
       

\bibitem{ir1} M. B. Aroux\'et, N. E. Echebest, and E. A. Pilotta,
   \textit{Inexact Restoration method for nonlinear optimization without
  derivatives,} J. Comput. Appl. Math. 290, pp. 26--43, 2015.

 \bibitem{ayvaz} M. T. Ayvaz, 
  \textit{A linked simulation-optimization model
  for simultaneously estimating the Manning's surface roughness values
  and their parameter structures in shallow water flows,} J. Hydrol. 500, pp. 183--199, 2013.
 


\bibitem{ir13} N. Banihashemi and C. Y. Kaya,
  \textit{ Inexact Restoration for
  Euler discretization of box-constrained optimal control problems, }J. Optim. Theory Appl. 156,
  pp. 726--760, 2013.

\bibitem{ir20} N. Banihashemi and C. Y. Kaya,  \textit{Inexact restoration and
  adaptive mesh refinement for optimal control, }J. Ind. Manag. Optim. 10, pp. 521--542, 2014.

\bibitem{inexact1} S. Bellavia, G. Gurioli, B. Morini, and
  Ph.L. Toint, \textit{Adaptive regularization algorithms with inexact
  evaluations for nonconvex optimization,}SIAM J. Optim. 29, pp. 2881--2915, 2019.
  
\bibitem{inexact2} S. Bellavia, G. Gurioli, B. Morini, and
  Ph.L. Toint, \textit{High-order Evaluation Complexity of a Stochastic
  Adaptive Regularization Algorithm for Nonconvex Optimization Using
  Inexact Function Evaluations and Randomly Perturbed Derivatives,}
  arXiv preprint arXiv:2005.04639, 2020.
 
\bibitem{ir23} S. Bellavia, N. Kreji\'c, and B. Morini,  \textit{Inexact
  restoration with subsampled trust‐region methods for finite‐sum
  minimization, }Comput. Optim. Appl.
  76, pp. 701--736, 2020.
      
 
 \bibitem{bkmr} S. Bellavia, N. Kreji\'c, B. Morini, and S. Rebegoldi,  
 \textit{A stochastic first-order trust-region method with inexact restoration for finite-sum minimization,} Comput. Optim. Appl. 84, pp. 53-84, 2023. 




\bibitem{bcrz} A. S. Berahas, F. E. Curtis, D. P. Robinson, and B.  Zhou,  {\it Sequential quadratic optimization for nonlinear equality 
constrained stochastic optimization,} SIAM J. Optim. 31, pp. 1352-1379, 2021.



 
\bibitem{ir15} E. G. Birgin, L. F. Bueno, and J. M. Mart\'{\i}nez,
   \textit{Assessing the reliability of general-purpose Inexact Restoration
  methods, }J. Comput. Appl. Math.
  282, pp. 1--16, 2015.
 
\bibitem{bkm2018} E. G. Birgin, N. Kreji\'c, and J. M. Mart\'{\i}nez,
   \textit{On the employment of Inexact Restoration for the minimization of
  functions whose evaluation is subject to errors, }Math. Comput. 87, pp. 1307--1326, 2018.
 
\bibitem{bkm2019} E. G. Birgin, N. Kreji\'c, and  J. M. Mart\'{\i}nez,
\textit{  Iteration and evaluation complexity on the minimization of functions
  whose computation is intrinsically inexact, }Math. Comput. 89, pp. 253--278, 2020.

  
\bibitem{bkm2021} E. G. Birgin, N. Kreji\'c, and J. M. Mart\'{\i}nez,  \textit{Economic inexact 
restoration for derivative-free expensive function minimization 
and applications,} J. Comput. Appl. Math. 410, 114193, 2022.

  

 
\bibitem{ir24} E. G. Birgin, R. D. Lobato, and J. M. Mart\'{\i}nez,
\textit{  Constrained optimization with integer and continuous variables using
  inexact restoration and projected gradients, }Bull. Comput. Appl. Math. 4, pp. 55--70, 2016.
 
\bibitem{ir9} E. G. Birgin and J. M. Mart\'{\i}nez, 
  \textit{Local convergence
  of an Inexact-Restoration method and numerical experiments,}J. Optim. Theory Appl. 127,
  pp. 229--247, 2005.



\bibitem{bfms} L. F. Bueno, A. Friedlander, J. M. Mart\'{\i}nez, and
  F. N. C. Sobral,\textit{ Inexact Restoration method for derivative-free
  optimization with smooth constraints, }SIAM J. Sci. Comput. 23, pp. 1189--1231, 2013.
 
\bibitem{ir16} L. F. Bueno, G. Haeser, and J. M. Mart\'{\i}nez, 
  \textit{A
  flexible Inexact-Restoration method for constrained optimization,}J. Optim. Theory. Appl. 165,
  pp. 188--208, 2015.

\bibitem{ir19} L. F. Bueno, G. Haeser, and J. M. Mart\'{\i}nez, 
  \textit{An
  inexact restoration approach to optimization problems with
  multiobjective constraints under weighted-sum scalarization,}Optim. Lett. 10, pp. 1315--1325, 2016.

\bibitem{ir22} L. F. Bueno and J. M. Mart\'{\i}nez,
  \textit{ On the complexity
  of an Inexact Restoration method for constrained optimization,}SIAM J. Optim. 30, pp. 80--101, 2020.



 \bibitem{cor} F. E. Curtis, M. J. O'Neill, and D. P. Robinson, \textit{Worst-case complexity of an SQP method for nonlinear equality constrained
 optimization}, COR@L Technical Report  21T-015, Lehigh University, Jan 6 2022.


\bibitem{crz}  F. E. Curtis, D. P. Robinson, and B. Zhou, \textit{Inexact Sequential Quadratic Optimization 
  for minimizing a stochastic objective function subject to deterministic nonlinear equality
   constraints, } 
  COR@L Technical Report  22T-01, Lehigh University, July 9  2021.

 
 
 
\bibitem{ir17} N. Echebest, M. L. Schuverdt, and R. P. Vignau, \textit{ An
  inexact restoration derivative-free filter method for nonlinear
  programming,}Comput. Appl. Math. 36,
  pp. 693--718, 2017.
 


\bibitem{ir18} P. S. Ferreira, E. W. Karas, M. Sachine, and 
  F. N. C. Sobral, 
  \textit{Global convergence of a derivative-free inexact
  restoration filter algorithm for nonlinear programming,}Optimization 66, pp. 271--292, 2017.
 
\bibitem{ir3} A. Fischer and A. Friedlander, 
  \textit{A new line search inexact
  restoration approach for nonlinear programming,}Comput. Optim. Appl. 46,
  pp. 336--346, 2010.

\bibitem{ir25} J. B. Francisco, D. S. Gon\c{c}alves, F. S. V. Baz\'an,
  and L. L. T. Paredes, \textit{Non-monotone inexact restoration method for
  nonlinear programming, }Comput. Optim. Appl. 76, pp. 867--888, 2020.

\bibitem{ir27} J. B. Francisco, D. S. Gon\c{c}alves, F. S. V. Baz\'an,
  and L. L. T. Paredes, \textit{Nonmonotone inexact restoration approach for
  minimization with orthogonality constraints, }Numer. Algorithms 86, pp. 1651--1684, 2021.
 
\bibitem{ir14} J. B. Francisco, J. M. Mart\'{\i}nez, L. Mart\'{\i}nez,
  and F. Pisnitchenko,
  \textit{ Inexact restoration method for minimization
  problems arising in electronic structure calculations, }Comput. Optim. Appl. 50,
  pp. 555--590, 2011.
 
\bibitem{ir12} M. A. Gomes-Ruggiero, J. M. Mart\'{\i}nez, and
  S. A. Santos, 
  \textit{Spectral Projected Gradient method with Inexact
  Restoration for minimization with nonconvex constraints,}SIAM J. Sci. Comput. 31, pp. 1628--1652,
  2009.


 \bibitem{gkv} C. C. Gonzaga, E. W. Karas, and M. Vanti, 
 
 \emph{A globally convergent filter method for nonlinear programming, }SIAM J. Optim.~14 pp.646-669, 2003.
 

\bibitem{inexact7} S. Gratton, E. Simon, D. Titley-Peloquin, and 
  Ph. L. Toint,  \textit{Minimizing convex quadratics with variable precision
  conjugate gradients,} Numer. Linear Algebra Appl. 28, e2337, 2021.
 
\bibitem{inexact4} S. Gratton, E. Simon, and Ph. L. Toint,
  \textit{ An
  algorithm for the minimization of nonsmooth nonconvex functions
  using inexact evaluations and its worst-case complexity, }Math. Program. 187, 1--24, 2021.
 
\bibitem{inexact6} S. Gratton and Ph. L. Toint,
  \textit{ A note on solving
  nonlinear optimization problems in variable precision,}Comput. Optim. Appl. 76,
  pp. 917--933, 2020.
  
\bibitem{ir6} E. Karas, E. A. Pilotta, and A. Ribeiro, 
  \textit{Numerical
  comparison of merit function with filter criterion in Inexact
  Restoration algorithms using hard-spheres problems,} Comput. Optim. Appl. 44,
  pp. 427--441, 2009.

\bibitem{ir10} C. Y. Kaya,  \textit{Inexact Restoration for Runge–Kutta
  discretization of optimal control problems,}SIAM J. Numer. Anal. 48, pp. 1492--1517, 2010.
 
\bibitem{ir8} C. Y. Kaya and J. M. Mart\'{\i}nez, \textit{Euler discretization
  and Inexact Restoration for optimal control, }J. Optim. Theory. Appl. 134, pp. 191--206, 2007.



 
\bibitem{inexact3} D. P. Kouri, M. Heinkenschloss, D. Ridzal, and
  B. G. van Bloemen Waanders,  \textit{Inexact objective function evaluations
  in a trust-region algorithm for PDE-constrained optimization under
  uncertainty,} SIAM J. Sci. Comput. 36,
  pp. A3011--A3029, 2014.

\bibitem{nkjmm} N. Kreji\'c and J. M. Mart\'{\i}nez,\textit{ Inexact
  Restoration approach for minimization with inexact evaluation of the
  objective function, }Math. Comput. 85,
  pp. 1775--1791, 2016.
 


\bibitem{leveque} R. J. LeVeque, {\it Finite Difference Methods for Ordinary and Partial Differential 
Equations. Steady-State and Time-Dependent Problems}, SIAM Publications, Philadelphia, 2007.



\bibitem{ir4} J. M. Mart\'{\i}nez, \textit{Inexact-restoration method with
  Lagrangian tangent decrease and new merit function for nonlinear
  programming, }J. Optim. Theory. Appl. 111, 39--58, 2001.
 
\bibitem{ir2} J.M. Mart\'{\i}nez and E. A. Pilotta, \textit{ Inexact
  restoration algorithms for constrained optimization, }J. Optim. Theory Appl. 104, pp. 135--163, 2000.

\bibitem{ir5} J. M. Mart\'{\i}nez and E. A. Pilotta,\textit{ Inexact
  restoration methods for nonlinear programming: advances and
  perspectives, }in L. Qi, K. Teo, X. Yang (eds.), \textit{Optimization
    and Control with Applications}, App. Optim. 96,
  Springer, New York, NY, 2005, pp. 271--291.

 \bibitem{mhh} A. Miele, H. Y. Huang, and J. C. Heideman,  {\it Sequential
gradient-restoration algorithm for the minimization of constrained
functions, ordinary and conjugate gradient version,}J. Optim. Theory Appl. 4, pp. 213-246, 1969.
 


\bibitem{livrocriab} J. L. Picanço, J. M. Mart\'{\i}nez, C. Pfeiffer, and J. F. Meyer (editors), {\it Conflitos,
  Riscos e Impactos Associados a Barragens}, CRIAB Publication, Institute of Advanced Studies of University of 
   Campinas, 2023.
 
\bibitem{rosen} J. B. Rosen,   {\it The gradient projection method for
nonlinear programming, Part 2, Nonlinear constraints,} SIAM J. on
Appl. Math.. 9, pp. 514-532, 1961.


\bibitem{saintvenant} A. J. C. Saint-Venant, \textit{Th\'eorie du mouvement
  non-permanent des eaux, avec application aux crues des rivi\`ere at
  \`a l'introduction des mar\'ees dans leur lit, }Comptes
    Rendus des S\'eances de Acad\'emie des Sciences 73, pp. 147--154,
  1871.
 

 \bibitem{sxxn} H-J M. Shi, Y. Xie, M. Q. Xuan, and J. Nocedal, 
 
{\it Adaptive Finite-Difference Interval Estimation for Noisy Derivative-Free Optimization, }SIAM J. Sci. Comput. 4, 2022. (doi: 10.1137/21M1452470)
 



 


\bibitem{ir21} C. E. P. Silva and M. T. T. Monteiro,  \textit{A filter
  inexact-restoration method for nonlinear programming,} TOP
  16, pp. 126--146, 2008.

\bibitem{ir26} J. Walpen, P. A. Lotito, E. M. Mancinelli, and
  L. Parente, \textit{ The demand adjustment problem via inexact restoration
  method,} Comput. Appl. Math. 39, article
  number 204, 2020.


\end{thebibliography}

\end{document}